\documentclass[nonblindrev]{informs3}

\OneAndAHalfSpacedXI

\usepackage{amsmath,amssymb,amsfonts}

\usepackage{natbib}
 \bibpunct[, ]{(}{)}{,}{a}{}{,}%
 \def\bibfont{\small}%

\usepackage{rotating}
\usepackage{fancyvrb}
\newtheorem{thm}{Theorem}

\newtheorem{rmrk}{Remark}
\usepackage{xcolor}

\definecolor{color1}{rgb}{0.95, 0.83, 0.31} 
\definecolor{color2}{rgb}{0.31, 0.71, 0.91} 
\definecolor{color3}{rgb}{0.47, 0.64, 0.26} 
\definecolor{color4}{rgb}{0.0, 0.35, 0.65} 
\definecolor{color5}{rgb}{0.17, 0.24, 0.56} 
\definecolor{color6}{rgb}{0.78, 0.13, 0.15} 
\definecolor{color7}{rgb}{0.0, 0.0, 0.0} 
\usepackage{algorithmic}            
\usepackage{algorithm}
\usepackage{multirow}

\usepackage{tikz,pgfplots}
\usetikzlibrary{matrix}
\usepgfplotslibrary{groupplots}
\pgfplotsset{compat=newest}
\usetikzlibrary{calc, positioning, shapes, fit, arrows, shadows, chains}
\tikzstyle{line} = [ draw, -latex']
\usetikzlibrary{shapes, arrows}
\usepgfplotslibrary{units}

\usepackage{booktabs} 

\usepackage{graphicx}  
\usepackage{subcaption} 

\usepackage{booktabs,caption}
\usepackage[flushleft]{threeparttable}

\TheoremsNumberedThrough     
\ECRepeatTheorems
\JOURNAL{Manufacturing \& Service Operations}

\EquationsNumberedThrough    

\MANUSCRIPTNO{IJDS-0001-1922.65}

\begin{document}

\RUNAUTHOR{}

\RUNTITLE{Hydrogen Network Expansion Planning considering the Chicken-and-egg Dilemma and Market Uncertainty}

\TITLE{Hydrogen Network Expansion Planning considering the Chicken-and-egg Dilemma and Market Uncertainty}


\ARTICLEAUTHORS{%
\AUTHOR{Sezen Ece Kayacik}
\AFF{Department of Operations, Faculty of Economics and Business, University of Groningen, the Netherlands,\EMAIL{s.e.kayacik@rug.nl}}
\AUTHOR{Beste Basciftci}
\AFF{Department of Business Analytics, Tippie College of Business, University of Iowa, Iowa City, Iowa, the United States, \EMAIL{beste-basciftci@uiowa.edu}}

\AUTHOR{Albert H. Schrotenboer}
\AFF{Operations, Planning, Accounting \& Control Group, School of Industrial Engineering, Eindhoven University of Technology, the Netherlands, \EMAIL{a.h.schrotenboer@tue.nl}}

\AUTHOR{Iris F. A. Vis, Evrim Ursavas}
\AFF{Department of Operations, Faculty of Economics and Business, University of Groningen, the Netherlands, \{i.f.a.vis@rug.nl, e.ursavas@rug.nl\}}
}

\ABSTRACT{%
 Green hydrogen is thought to be a game changer for reaching sustainability targets. However, the transition to a green hydrogen economy faces a critical challenge known as the `chicken-and-egg dilemma', wherein establishing a hydrogen supply network relies on demand, while demand only grows with reliable supply. In addition, as the hydrogen market is in the early stage, predicting demand distributions is challenging due to lack of data availability. This paper addresses these complex issues through a risk-averse framework with the introduction of a distributionally robust hydrogen network expansion planning problem under decision-dependent demand ambiguity. The problem optimizes location and production capacity decisions of the suppliers considering the moments of the stochastic hydrogen demand as a function of these investment decisions. 
To obtain tractable representations of this problem, we derive two different reformulations that consider continuous and discrete hydrogen demand support sets under different forms of decision dependencies. To efficiently solve the reformulations, we develop a tailored algorithm based on the column-and-constraint generation approach, and enhance the computational performance through solving the master problems to a relative optimality gap, decomposing the subproblems, and integrating pre-generated columns and constraints. Comparative experiments demonstrate the superiority of our algorithm over benchmark approaches including the classical column-and-constraint generation algorithm and Benders decomposition with significant speedups. To validate the effectiveness of our approach in decision-making, we investigate a real case study leveraging data from the ``Hydrogen Energy Applications in Valley Environments for Northern Netherlands (HEAVENN)" project.
The results reveal that considering the chicken-and-egg dilemma under uncertain hydrogen market conditions leads to earlier and more diverse investments, providing critical insights for policymakers based on the degree of decision dependency.
}


\maketitle

\section{Introduction}
The European Union aims to achieve net zero emission targets by 2050 in line with global climate mitigation goals \citep{eu2012}. In this context, green hydrogen, produced through electrolysis using renewable energy sources like wind or solar power, has emerged as a promising energy option with its 
zero emissions \citep{hosseini2016hydrogen,glenk2019economics}. As shown in the European Union hydrogen roadmap, various sectors such as mobility, industry, and heating can benefit from green hydrogen in mitigating carbon emissions \citep{undertaking2019hydrogen}. To make green hydrogen cost-competitive with fossil fuel alternatives, electrolyzer capacity expansion is necessary to achieve economies of scale. However, the transition to a hydrogen economy faces the critical `chicken-and-egg dilemma': the establishment of a hydrogen supply network depends on demand, but demand, in turn, only grows with a reliable supply.
The central debate for hydrogen policymakers in many countries is, therefore, when and how hydrogen networks should be developed so that the transition toward a future green-hydrogen economy can be made sustainably and cost-effectively.





Hydrogen network expansion planning is a long-term initiative that is crucial for sustaining growth and meeting increasing hydrogen demand across various sectors \citep{wappler2022building}. Ongoing projects, such as the European Hydrogen Backbone (EHB 2020), the Appalachian Hydrogen Hub (ARCH 2023), and the Western Green Energy Hub (WGEH 2024), illustrate the global momentum toward establishing future hydrogen networks. A central aspect of this establishment involves determining optimal locations and production capacities of electrolyzers to supply hydrogen to customers efficiently. Identification of customer demand is important as it is a driving factor in determining the locations and capacities of electrolyzers. However, accurately predicting hydrogen demand presents a significant challenge because of its inherent uncertainty \citep{iea2022global}. Furthermore, unlike established markets that rely on historical data for demand predictions, the early-stage hydrogen market lacks this historical context, making demand predictions even more challenging.


%


The challenge of demand prediction does not stand on itself, as it is embedded within another critical challenge, namely the  `chicken-and-egg dilemma' corresponding to the infrastructure-dependent adoption of emerging technologies as in energy transition \citep{Ma2019_Chickenegg, Brozynski2022_Chickenegg}. For the widespread adoption of hydrogen as a clean energy source, there needs to be a well-established production and distribution network to supply hydrogen. Nevertheless, investing in such network is often seen as economically viable only when there is sufficient demand for hydrogen. Paradoxically, the demand for hydrogen is likely to grow when there is a reliable supply \citep{worldenergy2022}. Navigating this dilemma poses a complex task for policymakers. Effective decision-making should recognize that the availability of hydrogen supply can drive the market to adopt hydrogen while simultaneously growing demand will cause businesses to expand supply capacities \citep{irenapolicy2022}. According to a report by \citet{iea2022global}, only 5\% of the announced hydrogen projects have taken firm investment decisions due to future demand uncertainties and the lack of available infrastructure. To overcome these challenges, mathematical analysis is essential for informing decision-making, encouraging projects to take firm steps and move forward with investment. Otherwise, the hydrogen market will face difficulties in scaling up.

To address these challenges, in this paper, we introduce a two-stage distributionally robust hydrogen network expansion planning problem, 
where the 
uncertain hydrogen demand depends on the investment decisions of the hydrogen production facilities, tackling the chicken-egg dilemma through the integration of decision-dependent uncertainties. 
The first stage of this problem involves determining optimal investment decisions, specifically the locations and production capacities of the hydrogen production facilities. Distributionally robustness ensures that hydrogen demand can come from a set of potential distributions whose mean is a function of these investment decisions, which is suitable for the inherent uncertainty of early-stage hydrogen market demand. Once demand is realized, the second stage determines operational decisions, including local production, transportation, and importation, to meet this demand. The goal is to minimize the expected sum of investment costs and operational costs under the worst-case demand distribution. To efficiently solve the resulting model, we develop a tailored algorithm based on column-and-constraint generation. We propose acceleration strategies to improve the performance of the algorithm. We show that our developed model and associated solution approach contribute to the central debate of policy-makers on how and when to develop hydrogen network. For this, we conduct a real case study leveraging data from  Europe’s first hydrogen valley, which is being built in the ``Hydrogen Energy Applications in Valley Environments for Northern Netherlands (HEAVENN)" project in the Netherlands.

In the remainder of this section, we discuss related work to highlight our theoretical and managerial contributions to the literature. We then summarize our contributions, outlining how our optimization framework addresses the challenges associated with the transition to a green hydrogen economy.




\subsection{Literature}





In real-world decision-making, uncertainties are prevalent and challenging to navigate. Two of the primary approaches addressing decision-making problems under  uncertainties differ in terms of how they capture these uncertainties. Specifically, stochastic programming \citep{wallace2005applications, birge2011introduction} aims to optimize the expected outcomes based on known probability distributions of uncertainties, whereas robust optimization \citep{ben2009robust,gorissen2015practical} takes a conservative approach by optimizing for the worst case within a predefined uncertainty set. Recognizing the inherent limitations of these approaches led to the development of a compromise known as distributionally robust optimization \citep{delage2010distributionally,rahimian2019distributionally}. This recent approach seeks optimal solutions for the worst-case probability distribution within a set of potential distributions, known as the ambiguity set.
In the literature, ambiguity sets can be defined in different ways depending on the level of information available. One common way is to use statistical-distance-based ambiguity sets that focus on distributions within a certain statistical distance of a reference distribution. The Wasserstein metric \citep{mohajerin2018data,kuhn2019wasserstein} and $\phi$-divergence \citep{ghosh2019robust} 
are example distance metrics for this category. These methods rely on historical data to establish a reference distribution. Since the hydrogen market is in the early development phase, only limited historical data is available, which hampers the creation of reliable reference distributions needed to employ these methods effectively.
Another common way is generalized-moment-based ambiguity sets \citep{basciftci2021distributionally,yu2022multistage}, which focus on moments and their proximate ranges. Employing generalized moments requires less data than is needed for constructing reliable reference distributions. Therefore, as applied in this study, a moment-based ambiguity set is more suitable than a statistical-distance-based ambiguity set for uncertain hydrogen demand.

While numerous papers explore distributionally robust optimization across different applications such as portfolio optimization \citep{delage2010distributionally,rujeerapaiboon2016robust}, supply chain management \citep{klabjan2013robust,basciftci2023resource} and machine learning \citep{duchi2019variance,chen2018robust}, 
the integration of decision-dependent uncertainty  within this framework remains limited. Decision-dependent uncertainty in this context refers to the case where the decisions can affect the underlying stochastic processes associated with the uncertain parameters. This concept is highly relevant in practical applications where system parameters, especially in new service industries or emerging markets like those for green hydrogen, exhibit uncertain features directly linked to prior decisions. For example, for fuel cell electric vehicles, customer demand varies significantly with the accessibility of charging stations, similar to the dependency of car-sharing services on the availability of service spots. Therefore, infrastructure significantly influences the demand for goods or services, but understanding the precise mechanisms is challenging due to the limited availability of historical data. Incorporating decision-dependent uncertainty into distributionally robust optimization could offer a promising solution for tackling this challenge.

Several papers consider decision-dependent uncertainties within distributionally robust optimization. \cite{luo2020distributionally} study five different ambiguity sets, including moment-based and distance-based variants. Under discrete support sets, they provide numerical examples on a small-scale instance for a newsvendor problem. \cite{basciftci2021distributionally} study distributionally robust optimization for a two-stage facility location problem under demand ambiguity, where demand is dependent on the first-stage location decisions. They provide a reformulation for a discrete support set of demand. \cite{noyan2022distributionally} focus on decision-dependent ball-centered ambiguity sets, integrating total variation distance and Wasserstein metrics. They provide alternative reformulations for the resulting nonlinear model to apply on machine scheduling and humanitarian logistics. Recently, \cite{YuBasciftci2024_MultimodalDDDR} integrate decision-dependent multimodal ambiguities to this problem by considering different modes of distributions with an additional layer of uncertainty around mode probabilities under moment-based and distance-based ambiguity sets.  \cite{doan2022distributionally} propose a reformulation for  a distributionally robust decision-dependent retrofitting problem where the survival probability depends on retrofitting decisions. They develop a constraint generation algorithm tailored to the retrofitting problem. \cite{yu2022multistage} provide reformulations for three different ambiguity sets for general multistage distributionally robust optimization. They solve the resulting reformulations by applying stochastic dual dynamic programming.
In the majority of the aforementioned studies, which consider discrete support sets, the size of support sets is often kept small to maintain computational tractability. However, this approach can sometimes limit the accuracy of uncertainty estimation due to the restricted size of the sets.
In addition, continuous support sets often lead to semi-infinite optimization problems that includes infinitely many constraints. Addressing this challenge requires developing solution algorithms as conventional optimization solvers cannot directly solve such formulations. We contribute to the literature by developing a solution algorithm to tackle problems with both discrete and continuous support sets.






For solving large-scale distributionally robust optimization problems, separation-based methods such as cutting plane-based algorithms \citep{ardestani2018value} 
 and column-and-constraint generation (C\&CG) algorithms \citep{zeng2013solving} can be used. Cutting-plane-based algorithms work by iteratively adding constraints to narrow down the feasible region, while the C\&CG algorithm iteratively generates decision variables and constraints by focusing on a subset of scenarios crucial for achieving an optimal solution robust against various uncertainties. 
To improve its computational performance, \cite{shehadeh2021distributionally,agra2018robust} propose problem-specific enhancements, whereas \cite{tsang2023inexact} introduce an inexact C\&CG method in which the master problem is solved to a relative optimality gap. Papers studying distributionally robust optimization with decision-dependent uncertainty typically propose reformulations without accompanying solution algorithms or direct application of existing solution techniques by mainly considering discrete support sets. 
Given that decision-dependent uncertainty introduces additional complexity, there is a need for enhanced algorithms capable of effectively solving large-scale instances of these complex problems. In this paper, we address this need by developing an advanced solution algorithm based on the C\&CG algorithm.

Effectively managing uncertainties is a core challenge in energy system problems like network expansion planning. We review the relevant literature that develops various optimization methodologies to address these uncertainties and support decision-making in network expansion planning. One line of research utilizes stochastic programming approaches to manage these uncertainties. For example, \cite{zou2018partially} propose partially adaptive stochastic optimization and \cite{Chen2024_ESEP} present multistage stochastic programming with chance-constraints for power generation expansion planning to determine optimal construction decision of generators. Another study by \cite{kayacik2024towards} propose a two-stage stochastic optimization model to determine optimal location and capacities for renewable energy sources and hydrogen storage. Another line of research focuses on robust optimization. For instance, \cite{moret2020decision} study a robust optimization framework to support decisions in strategic energy planning considering various uncertainties and \cite{abdin2022optimizing} study a long-term generation expansion planning problem using a multistage adaptive robust approach. Recently, distributionally robust optimization approaches have been integrated into these problem settings, offering a new perspective on managing uncertainties. \cite{pourahmadi2019distributionally} propose a distributionally robust chance-constrained optimization problem for a centralized generation expansion planning problem and \cite{guevara2020machine} develop a combined distributionally robust optimization and machine learning approach to determine investment decisions in energy planning problems. A recent study by \cite{qiu2023decision} proposes a distributionally robust integrated generation, transmission, and storage expansion problem considering a decision-dependent distance-based ambiguity set to reflect the effect of the sizing of renewable energy sources on their energy output. Despite these advancements, the underlying decision-dependent uncertainties that arise as part of the chicken-egg problem along with an efficient algorithm for its solution and its practical implications have not been addressed in the relevant literature.


\subsection{Contributions}
 We summarize the contributions of this paper as follows:
\begin{itemize}
    \item We introduce a novel two-stage distributionally robust decision-dependent hydrogen network expansion planning problem, where the mean of uncertain hydrogen demand depends on investment decisions. In the first stage, the model determines the location and production capacity decisions for the hydrogen facilities. Once the decision-dependent demand is realized, the model then determines operational decisions such as local production, import, and transportation in the second stage.  
    \item We provide a monolithic reformulation for the proposed problem under a generic decision-dependent moment function, which can be used for continuous and discrete support sets. Then, we provide two types of piecewise linear moment functions, depending on the location and 
    production capacity decisions for capturing the endogenous uncertainty. We derive additional reformulations specific to these moment functions.
    \item We develop a solution approach based on a column-and-row generation algorithm, and propose various computational enhancements. Our acceleration strategies include solving the master problem to a relative optimality gap, decomposing  the subproblem by exploiting the problem structure, and including pre-generated columns and constraints. 
    In our computational experiments, we demonstrate the superiority of our algorithm across various instances, successfully solving 75\% of them to optimality within a 5-hour time limit with significant speedups. In contrast, for the benchmark solution approaches, the classical C\&CG algorithm solves 41\% of the instances, while the Benders decomposition algorithm solves only  35\%.
    \item We validate our approach by means of a case study based on real-world data to create hydrogen network expansion plans for the Northern Netherlands by providing insights on the i) value of decision dependency over investment plans, ii) out-of-sample performances of the proposed approach against benchmark optimization approaches, iii) sensitivity analyses on problem parameters, and iv) integration of location-based and capacity-based moment functions for capturing the decision dependency. We show that considering the chicken-and-egg dilemma leads to earlier and more diverse investment decisions, increasing the net profit of the hydrogen network and accelerating the adoption of green hydrogen.
\end{itemize}

The remainder of the paper is organized as follows. Section~\ref{sec:Model} introduces the distributionally robust decision-dependent optimization problem and provides its  reformulations.  Section~\ref{sec:SolutionApproach} outlines our solution approach based on C\&CG algorithm. Section~\ref{sec:ComputationalExperiments} provides the computational and case study results, demonstrating the outperforming performance of our solution approach and offering managerial insights from the real case study. Lastly, Section~\ref{sec:Conclusions} concludes the paper with our final remarks.

\section{Problem Description} \label{sec:Model}

In this section, we introduce a distributionally robust hydrogen network expansion planning problem with decision-dependent uncertainty that incorporates the chicken-and-egg dilemma and uncertain hydrogen market conditions. The problem involves designing a hydrogen supply network consisting of local production facilities and ports, where hydrogen is either produced locally or imported via ports to meet demand at various customer locations. The goal is to establish a long-term hydrogen network expansion plan. In the first-stage, we determine optimal locations and production capacities of hydrogen facilities for each period over a finite planning horizon. Once these facilities are established, hydrogen demand realizes in the second-stage from a set of potential distributions whose mean is a function of the first-stage location and capacity decisions. After the realization of this \textit{decision-dependent demand}, the amount of local production and importation are determined to transport hydrogen and meet demand across various locations. The goal is to minimize the distributionally worst-case expected sum of investment and operational costs.



The rest of this section is structured as follows: Section~\ref{Model: formulation} outlines the problem description and introduces a distributionally robust optimization model with the ambiguity set. Subsequently, Section~\ref{Model: moments} discusses different decision-dependent moment functions. Lastly, Section~\ref{Model: reformulations} provides two reformulations of the distributionally robust optimization model under discrete and continuous support sets of decision-dependent hydrogen demand.

\subsection{Model Formulation}\label{Model: formulation}

We consider the expansion planning of a hydrogen network over periods $t \in \mathcal{T}$. The network nodes $i \in \mathcal{N}$ are partitioned into supply nodes to open a hydrogen facility ($\mathcal{N}^S \subseteq \mathcal{N}$), port nodes for import of hydrogen ($\mathcal{N}^I \subseteq \mathcal{N}$), and demand nodes ($\mathcal{N}^D \subseteq \mathcal{N}$). We define binary variables  $x_{it}$  equaling 1 if a hydrogen facility is open at a candidate node $i \in \mathcal{N}^S$ at time $t \in \mathcal{T}$. The variable $y_{it}$ identifies the production capacity of a hydrogen facility at a location $i \in \mathcal{N}^S$ at time $t \in \mathcal{T}$. The production capacity  of a hydrogen facility at location  $i \in \mathcal{N}^S$  at time $t \in \mathcal{T}$ can be limited by the capacity limit $ \bar C_{it}$. We incur a one-time set-up cost $c^{F}_{it}$ and a variable production capacity cost $c^{V}_{it}$ for the investment of hydrogen facilities.

After the opening and production capacity decisions take place, random hydrogen demand $\xi_{jt}$ for location $j \in \mathcal{N}^D$ at time  $t \in \mathcal{T}$ is realized. 
After demand realizes, we determine the amount of hydrogen  $h_{it}$ produced at node $i \in \mathcal{N}^S $ at time $t \in \mathcal{T}$ and the amount of hydrogen $v_{it}$ imported from port  $i \in \mathcal{N}^I $ at time $t \in \mathcal{T}$.
To directly transport the produced or imported hydrogen to demand nodes, we determine the amount of hydrogen $z_{ijt}$  transported between nodes $i,j \in \mathcal{N} $ at time $t \in \mathcal{T}$. We incur production costs $c^{P}_{it}$ for per kilogram of hydrogen produced at supply nodes $i \in \mathcal{N}^S$ at time $t \in \mathcal{T}$ and import costs $c^{I}_{it}$ for per kilogram of hydrogen imported from port nodes $i \in \mathcal{N}^I$ at time $t \in \mathcal{T}$. Based on the distance ($\Delta_{ij}$) between two nodes $i,j \in \mathcal{N}$, we incur transportation cost  $c^T_{ijt}$ for per kilogram of hydrogen at time $t \in \mathcal{T}$.
We assign a revenue $R_{jt}$ for selling per kilogram of hydrogen at node $j \in \mathcal{N}^D$ at time $t \in \mathcal{T}$. We note that, throughout the paper, matrix and vector notations are denoted in bold. 

Due to the early phase of the hydrogen market, the actual demand distribution is unknown. To this end, we assume that we only know the set of possible demand distributions, whose means $\mu_{jt}(\mathbf{x,y})$ depend on the investment decisions. Accordingly, we  define a support set $\mathcal{S} = \{  \underline{\xi}_{jt} \leq \xi_{jt} \leq \overline{\xi}_{jt}, \forall j \in \mathcal{N}^D, t \in \mathcal{T} \}$, where $\underline{\xi}_{jt}$ and $\overline{\xi}_{jt}$ correspond to lower and upper bounds of the corresponding demand values, respectively. Each element, denoted as $\mathbb{P}$, in the probability set $P(\mathcal{S})$, represents a probability measure associated with $\boldsymbol{\xi}$. Then, we construct the ambiguity set $U(\mathbf{x,y})$ depending on the investment decisions $\mathbf{x,y}$ as follows:
\begin{equation}\label{eq:U(x)}
 U(\mathbf{x,y}):=\left\{\mathbb{P} \in \mathcal{P}(\mathcal{S}): \begin{array}{l}
\mathbb{P}(\boldsymbol{\xi} \in \mathcal{S})=1, \\
\mu_{jt}(\mathbf{x,y}) - \epsilon_{jt} \leq \mathbf{E}_{\mathbb{P}}\left[\xi_{jt}\right] \leq \mu_{jt}(\mathbf{x,y}) + \epsilon_{jt}, \ \ \forall j \in \mathcal{N}^D, t \in \mathcal{T}.
\end{array}\right\},
\end{equation}
where the first constraint ensures that the sum of probabilities over the support set equals 1, while the second constraint controls the mean ensuring that the true mean is within an $\ell_1$-based distance $\epsilon_{jt}$ to the mean $\mu_{jt}(\mathbf{x,y})$. If the mean is perfectly known, we can set $\boldsymbol{\epsilon}$ to zero. Note that moment function $\mu_{jt}(\mathbf{x,y})$ is presented in its generic form in the ambiguity set \eqref{eq:U(x)}. However, in Section \ref{Model: moments}, we provide specific moment functions depending on the investment decisions to present insights on their construction. 

Then, we formulate the distributionally robust decision-dependent hydrogen network expansion problem as follows:
\begin{subequations} \label{RobustModel}
\begin{align}
    \min_{\mathbf{x},\mathbf{y}} & \quad \sum_{t \in \mathcal{T}} \sum_{i \in \mathcal{N}^S} \left( c^{F}_{it} (x_{it}- x_{i(t-1)} )+ c^{V}_{it} y_{it} \right) +\max _{\mathbb{P} \in U(\mathbf{x,y})} \mathbf{E}_{\mathbb{P}}[f(\mathbf{x},\mathbf{y},\boldsymbol{\xi})], \span  \span \label{RobustModel_obj} \\
    \text{s.t.} & \quad  y_{it}  \leq  \bar C_{it} x_{it},  \quad & \forall& i \in \mathcal{N}^S, t \in \mathcal{T}, \label{RobustModel_C1}  \\ 
     & \quad x_{it}  \geq  x_{i(t-1)},  \quad & \forall& i \in \mathcal{N}^S , t \in \mathcal{T} \backslash \{1\}, \label{RobustModel_C2} \\
     & \quad x_{it} \in \{0,1\}, y_{it} \geq 0,  \quad & \forall& i \in \mathcal{N}^S, t \in \mathcal{T}.  \label{RobustModel_C3}
\end{align}
\end{subequations}

Here, objective function~\eqref{RobustModel_obj} minimizes the set-up and capacity cost of hydrogen facilities, and the worst case expectation of the operational cost $f(\mathbf{x},\mathbf{y},\boldsymbol{\xi})$ considering the set of distributions from the ambiguity set \eqref{eq:U(x)}. Constraint~\eqref{RobustModel_C1} links the opening and capacity decisions for the set of supply nodes. Constraint~\eqref{RobustModel_C2} guarantees that any facility that is opened remains operational and cannot be closed until the end of the planning horizon. 

We formulate the operational second-stage problem $f(\mathbf{x},\mathbf{y},\boldsymbol{\xi})$ under given first-stage decisions $\mathbf{x},\mathbf{y}$ and demand realization $\boldsymbol{\xi}$ as follows:
\begin{subequations}\label{SubProblem}
\begin{align}
f(\mathbf{x},\mathbf{y},\boldsymbol{\xi})=\min_{\mathbf{h},\mathbf{z},\mathbf{v}} & \sum_{t \in \mathcal{T}} \left( \sum_{i \in \mathcal{N}^S} c^{P}_{it}h_{it}  + \sum_{i \in \mathcal{N}^I}  c^{I}_{it} v_{it} + \sum_{j \in \mathcal{N}^D} \left( \sum_{i \in \mathcal{N}^S \cup \mathcal{N}^I  } c^{T}_{ijt}   z_{ijt} -  R_{jt} \xi_{jt} \right) \right), \span \span \label{SP:obj}  \\ 
\text{s.t.} & \quad   h_{it}  \leq \sum_{t' \in [1,t]} y_{it'},  \hspace{3cm}
&\forall& i \in \mathcal{N}^S, t \in \mathcal{T}, \label{SP:C1}  \\ 
& \quad h_{it}  =  \sum_{j \in \mathcal{N}^D} z_{ijt}, \quad &\forall& i \in \mathcal{N}^S, t \in \mathcal{T},  \label{SP:C2}  \\ 
&  \quad  v_{it}  =  \sum_{j \in \mathcal{N}^D} z_{ijt}, \quad &\forall& i \in \mathcal{N}^I, t \in \mathcal{T}, \label{SP:C3}  \\ 
& \quad   \sum_{i \in \mathcal{N}^S \cup \mathcal{N}^I  } z_{ijt} = \xi_{jt}, \quad &\forall& j \in \mathcal{N}^D, t \in \mathcal{T},  \label{SP:C4} \\ 
& \quad \mathbf{h}, \mathbf{z}, \mathbf{v} \geq \mathbf{0}. & \label{SP:C5}
\end{align}\end{subequations}
Here, objective function~\eqref{SP:obj} minimizes the overall cost of producing, importing, and transporting hydrogen while considering the reward associated with meeting hydrogen demand. Constraint~\eqref{SP:C1} ensures that hydrogen production at each time period cannot be more than the capacity built until that period. Constraints~\eqref{SP:C2} and \eqref{SP:C3} control the hydrogen amounts distributed from supply and import nodes to a demand node with respect to production and import amounts, respectively. Constraint~\eqref{SP:C4} guarantees that hydrogen demand is satisfied through local production and import.  

Model~\eqref{RobustModel} includes the inner maximization problem $\max_{\mathbb{P} \in U(\mathbf{x,y})} \mathbf{E}_{\mathbb{P}}[f(\mathbf{x},\mathbf{y},\boldsymbol{\xi})]$, which can be formulated for  given $\mathbf{x},\mathbf{y}$ as follows:
\begin{subequations} \label{InnerModel}
    \begin{align}
\max_{\mathbb{P} \in \mathcal{P}(\mathcal{S})} & \quad  \int_{\mathcal{S}} f(\mathbf{x},\mathbf{y},\boldsymbol{\xi}) \mathbb{P}(d \boldsymbol{\xi}), \\
\text{s.t.} & \quad \int_{\mathcal{S}} \mathbb{P}(d \boldsymbol{\xi})=1 , \label{Inner:C1} \\
& \quad  \int_{\mathcal{S}} \xi_{jt} \mathbb{P}(d \boldsymbol{\xi}) \leq \mu_{jt}(\mathbf{x,y})+\epsilon_{jt},  & \forall j \in \mathcal{N}^D, t \in \mathcal{T}, \label{Inner:C2} \\
& \quad \int_{\mathcal{S}}  \xi_{jt} \mathbb{P}(d \boldsymbol{\xi}) \geq \mu_{jt}(\mathbf{x,y})-\epsilon_{jt},  & \forall  j \in \mathcal{N}^D, t \in \mathcal{T}.  \label{Inner:C3}
\end{align}\end{subequations}

The above reformulation of the inner maximization problem is obtained by using the definition of the ambiguity set $U(\mathbf{x,y})$ in Equation \eqref{eq:U(x)}. Constraint \eqref{Inner:C1}  ensures that the sum of probabilities across the support set equals 1, while Constraints \eqref{Inner:C2} and  \eqref{Inner:C3} control the mean.

To obtain a monolithic reformulation of Model \eqref{RobustModel}, we 
convert Model \eqref{InnerModel} into a minimization problem through its dual formulation. We introduce dual variables $\alpha, \beta^1_{jt}, \beta^2_{jt} $ for all $j \in \mathcal{N}^D$, $t \in \mathcal{T}$ corresponding to Constraints~\eqref{Inner:C1}, \eqref{Inner:C2}, and \eqref{Inner:C3}, respectively. Then, we derive the dual model as follows:
\begin{subequations} \label{DualInnerModel}
\begin{align}
\min_{\alpha, \boldsymbol{\beta}^1,\boldsymbol{\beta}^2} & \quad  \alpha + \sum_{t \in \mathcal{T}} \sum_{j \in \mathcal{N}^D}  \left(  \beta^1_{jt} (\mu_{jt}(\mathbf{x,y})+\epsilon_{jt}) -  \beta^2_{jt} (\mu_{jt}(\mathbf{x,y})-\epsilon_{jt}) \right), \span \\
\text{s.t.} & \quad   \alpha + \sum_{t \in \mathcal{T}} \sum_{j \in \mathcal{N}^D} (\beta^1_{jt} -\beta^2_{jt}) \xi_{it}  \geq f(\mathbf{x},\mathbf{y},\boldsymbol{\xi}),   & \forall \boldsymbol{\xi} \in \mathcal{S}, \label{DI: C1} \\
& \quad \boldsymbol{\beta}^1, \boldsymbol{\beta}^2 \geq 0. & \label{DI: C2}
\end{align}\end{subequations}

By replacing the inner maximization problem $\max_{\mathbb{P} \in U(\mathbf{x,y})} \mathbf{E}_{\mathbb{P}}[f(\mathbf{x},\mathbf{y},\boldsymbol{\xi})]$ in Model~\eqref{RobustModel} with its dual formulation \eqref{DualInnerModel}, an equivalent formulation of Model~\eqref{RobustModel} can be obtained as illustrated 
in the following theorem.
\begin{thm}  An equivalent formulation to Model~\eqref{RobustModel} can be obtained as follows:
\begin{align}  
\min_{\mathbf{x},\mathbf{y},\alpha,\boldsymbol{\beta}^1,\boldsymbol{\beta}^2} & \quad \sum_{t \in \mathcal{T}} \sum_{i \in \mathcal{N}^S} \left( c^{F}_{it} (x_{it}- x_{i(t-1)} )+ c^{V}_{it} y_{it} \right) +  \alpha + \sum_{t \in \mathcal{T}} \sum_{j \in \mathcal{N}^D}  \left(  \beta^1_{jt} (\mu_{jt}(\mathbf{x,y})+\epsilon_{jt}) -  \beta^2_{jt} (\mu_{jt}(\mathbf{x,y})-\epsilon_{jt}) \right), \label{EquivalentRobustModel}\\
    \text{s.t.} & \quad \eqref{RobustModel_C1}-\eqref{RobustModel_C3}, \eqref{DI: C1}-\eqref{DI: C2}.  \notag
\end{align}
\end{thm}

\subsection{Moment Functions} \label{Model: moments}
As discussed earlier, investments in new hydrogen facilities can influence customer demand by enhancing accessibility to hydrogen supply, thereby increasing customer interest in adopting hydrogen technology. To capture this effect, we first consider location-based decision-dependent moment functions, where only the location of supply affects demand, which is captured by a moment-function that only depends on $\mathbf{x}$. Secondly, we extend the decision dependency by considering the potential relationship between capacity decisions and demand, effectively having a moment function that depends only on $\mathbf{y}$.

\subsubsection{Location-based Moment Function.}

We consider that the demand at location $j$ at time $t$ increases from a base demand estimate $\bar \mu_{jt}$ when new supply facilities are opened in location $j$’s neighborhood \citep{basciftci2021distributionally}. The impact of opening location $i$ on the mean demand at location $j$ is controlled by the parameters $\lambda^{l}_{ij} \in [0,1]$ based on the distance between locations $i,j$. We let nearby locations have a higher impact on the mean, and further locations have less effect. Therefore, $\lambda^{l}_{ij}$ can be considered as a decreasing function of the distance between locations $i$ and $j$. We also consider a maximum potential demand increase limit  $\overline{B}_{jt}$ for each location and time period. Once this threshold is reached, the opening of new locations cannot further increase the demand. Then, we specify the location-based piecewise linear moment function $\mu_{jt}(\mathbf{x})$ for each period $t \in \mathcal{T}$ and location $j \in \mathcal{N}^D$ as follows:
\begin{equation} \label{Decision-dependent bounded moment}
\mu_{jt}(\mathbf{x})= \min \left( \bar \mu_{jt} \left(1 + \sum_{i \in \mathcal{N}^S} \lambda^{l}_{ij} x_{it} \right), \bar \mu_{jt} \left(1 + \overline{B}_{jt} \right)  \right). 
\end{equation}

\subsubsection{Capacity-based Moment Function.}
In addition to considering the effect of a supply location on customer demand, we introduce the setting in which the actual supply capacity affects hydrogen demand. 
This setting can arise when hydrogen demand is expected to be impacted by changes in the supply capacity levels. 
To formulate this setting, we define a set of capacity ranges $\mathcal{A}_{it}$ for each location $i \in \mathcal{N}^S$ and time $t \in \mathcal{T}$, where each range $r \in \mathcal{A}_{it}$  
has a lower $(l_{itr})$ and an upper bound $(u_{itr})$. 
Let $e_{itr}$ be a binary variable that equals to 1 if capacity $y_{it}$ is within $[l_{itr},u_{itr}]$, and 0 otherwise. Accordingly,  we modify Model~\eqref{RobustModel} by including constraints to identify the range of capacity at location $i \in \mathcal{N}^S$ and time $t \in \mathcal{T}$  as follows: 
\begin{subequations} \label{Capacity based moment constraints}
\begin{align}
\min_{\mathbf{x},\mathbf{y}} & \quad \sum_{t \in \mathcal{T}} \sum_{i \in \mathcal{N}^S} \left( c^{F}_{it} (x_{it}- x_{i(t-1)} )+ c^{V}_{it} y_{it} \right) +\max _{\mathbb{P} \in U(\mathbf{e})} \mathbf{E}_{\mathbb{P}}[f(\mathbf{x},\mathbf{y},\boldsymbol{\xi})], \span \span \\
\text{s.t.} & \quad \eqref{RobustModel_C1}-\eqref{RobustModel_C3},  \notag \\ 
& \quad l_{itr} - M(1-e_{itr}) \leq \sum_{t' \in [1,t]}  y_{it'} \leq  u_{itr} + M(1-e_{itr}),   & \forall&  r \in \mathcal{A}_{it},  i \in \mathcal{N}^S, t \in \mathcal{T}, \label{CapDDU: C1} \\
& \quad \sum_{r \in \mathcal{A}_{it}} e_{itr} = 1,  & \forall&    i \in \mathcal{N}^S, t \in \mathcal{T}, \label{CapDDU: C2} \\
& \quad e_{itr} \in \{0,1\}, & \forall&  r \in \mathcal{A}_{it}, i \in \mathcal{N}^S, t \in \mathcal{T}. \label{CapDDU: C3} 
\end{align}
\end{subequations}

Here, Constraint~\eqref{CapDDU: C1} identifies the range of the cumulative capacity of a hydrogen facility at location $i \in \mathcal{N}^S$ and time $t \in \mathcal{T}$. Constraint~\eqref{CapDDU: C2} ensures that the cumulative capacity can only belong to one range, thereby preventing simultaneous allocation to multiple ranges.

For the capacity-based scenario, we define the parameter $\lambda^{c}_{ijr} \in [0,1]$ for controlling the impact of capacity decisions on mean demand, which is a decreasing factor based on the distance between locations $i$ and $j$. Furthermore, as the capacity range $r$ increases, $\lambda^{c}_{ijr}$ also increases, indicating that higher capacity ranges have a greater impact on demand. Then, we define the capacity-based piecewise linear moment function $\mu_{jt}(\mathbf{e})$ for each period $t \in \mathcal{T}$ and location $j \in \mathcal{N}^D$ as follows:
\begin{equation} \label{Decision-dependent capacity moment}
\mu_{jt}(\mathbf{e})= \min \left( \bar \mu_{jt} \left(1 + \sum_{i \in \mathcal{N}^S} \sum_{r \in \mathcal{A}_{it}} 
 \lambda^{c}_{ijr} e_{itr} \right) , \bar \mu_{jt} \left(1 + \overline{B}_{jt} \right)\right).
\end{equation}

\subsection{Reformulations under Decision-dependent Moment Functions}\label{Model: reformulations}
To simplify notation throughout the remainder of this paper, we consider the unbounded version of the location-based moment function, where $ \boldsymbol{\overline{B}} = \boldsymbol{\infty}$. Thus, we derive our model reformulations based on Moment Function~\eqref{Decision-dependent moment}.  The reformulations for the bounded location-based moment function~\eqref{Decision-dependent bounded moment} and capacity-based moment function~\eqref{Decision-dependent capacity moment} are provided in Appendices~\ref{Appendix: Bounded} and \ref{Appendix: capacity-based}, respectively.
\begin{equation} \label{Decision-dependent moment}
\mu_{jt}(\mathbf{x})=  \bar \mu_{jt} \left(1 + \sum_{i \in \mathcal{N}^S} \lambda^{l}_{ij} x_{it} \right).
\end{equation}

We replace location-based moment function~\eqref{Decision-dependent moment} into Model~\eqref{EquivalentRobustModel} as follows:
\begin{align}
\min_{\mathbf{x},\mathbf{y}, \alpha,\boldsymbol{\beta}^1, \boldsymbol{\beta}^2} & \quad \sum_{t \in \mathcal{T}} \sum_{i \in \mathcal{N}^S} \left( c^{F}_{it} (x_{it}- x_{i(t-1)} )+ c^{V}_{it} y_{it} \right) +  \alpha + \sum_{t \in \mathcal{T}} \sum_{j \in \mathcal{N}^D}  \left(  \bar \mu_{jt} \beta^1_{jt} + (\bar \mu_{jt} \sum_{i \in \mathcal{N}^S} \lambda^{l}_{ij} x_{it}) \beta^1_{jt} \right. \nonumber \\
& \hspace{6cm} \left. + \epsilon_{jt} \beta^1_{jt} -  \bar \mu_{jt} \beta^2_{jt} - ( \bar \mu_{jt} \sum_{i \in \mathcal{N}^S} \lambda^{l}_{ij} x_{it}) \beta^2_{jt}  + \epsilon_{jt}\beta^2_{jt}) \right), \label{MomentReformulation} \\
\text{s.t.} & \quad \eqref{RobustModel_C1}-\eqref{RobustModel_C3}, \eqref{DI: C1}-\eqref{DI: C2}.  \nonumber 
\end{align}
The objective function of Model~\eqref{MomentReformulation} includes the bilinear terms $x_{it} \beta^1_{jt}$ and $x_{it} \beta^2_{jt}$. For linearization, we introduce auxiliary variables $\Phi^1_{ijt}$ and $\Phi^2_{ijt}$ for all $i \in \mathcal{N}^S, j \in \mathcal{N}^D, t \in \mathcal{T}$ and replace the bilinear terms with these auxiliary variables. Let $\bar \beta^m_{jt}$ be an upper bound on the variable $\beta^m_{jt}$ where $m \in \{1,2\}$. Consequently, we obtain the following McCormick constraints: 
\begin{subequations}\label{McCormicks}
\begin{align}
\Phi^m_{ijt} &\leq \bar \beta^m_{jt} x_{it},  & \forall m \in \{1,2\}, j \in \mathcal{N}^D, i \in \mathcal{N}^S, t \in \mathcal{T}, \label{MC1} \\
\Phi^m_{ijt} &\geq 0, & \forall m \in \{1,2\}, j \in \mathcal{N}^D, i \in \mathcal{N}^S, t \in \mathcal{T}, \label{MC2} \\
\Phi^m_{ijt} &\leq \beta^l_{jt}, & \forall m \in \{1,2\}, j \in \mathcal{N}^D, i \in \mathcal{N}^S, t \in \mathcal{T}, \label{MC3} \\
\Phi^m_{ijt} &\geq \beta^m_{jt} - \bar \beta^m_{jt} \left(1-x_{it}\right), \quad & \forall m \in \{1,2\}, j \in \mathcal{N}^D, i \in \mathcal{N}^S, t \in \mathcal{T}. \label{MC4}
\end{align}
\end{subequations}

Model~\eqref{MomentReformulation} is not computationally tractable in its present form since Constraint~\eqref{DI: C1} is infinitely many. Therefore, in addition to the inclusion of McCormick constraints, we update Constraint~\eqref{DI: C1} by taking the maximum of the right-hand side, resulting in the following equivalent formulation:
\begin{subequations} \label{ReformulationFirstStep}
\begin{align}
\min_{\mathbf{x},\mathbf{y},\alpha,\boldsymbol{\beta}^1, \boldsymbol{\beta}^2,\boldsymbol{\Phi}} & \quad  \sum_{t \in T} \sum_{i \in \mathcal{N}^S} \left( c^{F}_{it} (x_{it}- x_{i(t-1)} )+ c^{V}_{it} y_{it} \right)  + \alpha + \sum_{j \in \mathcal{N}^D} \sum_{t \in \mathcal{T}} \left(  \bar \mu_{jt} \beta^1_{jt} + \bar \mu_{jt} \sum_{i \in \mathcal{N}^S} \lambda^{l}_{ij} \Phi^1_{ijt} \right. \nonumber \\
& \hspace{6cm} \left. + \epsilon_{jt} \beta^1_{jt} -  \bar \mu_{jt} \beta^2_{jt} - \bar \mu_{jt} \sum_{i \in \mathcal{N}^S} \lambda^{l}_{ij} \Phi^2_{ijt}  + \epsilon_{jt}\beta^2_{jt} \right), \\
\text{s.t.} & \quad \eqref{RobustModel_C1}-\eqref{RobustModel_C3}, \eqref{DI: C2}, \eqref{MC1}-\eqref{MC4},  \nonumber \\
& \quad \alpha \geq \max_{\boldsymbol{\xi}\in \mathcal{S}} \left\{ f(\mathbf{x},\mathbf{y},\boldsymbol{\xi}) -  \sum_{j \in \mathcal{N}^D} \sum_{t \in \mathcal{T}}  (\beta^1_{jt} -\beta^2_{jt}) \xi_{jt} \right\}.\label{Constraint continuous}
\end{align}\end{subequations}

To eliminate the inner minimization in Constraint~\eqref{Constraint continuous} within $f(\mathbf{x},\mathbf{y},\boldsymbol{\xi})$, we take the dual of  $f(\mathbf{x},\mathbf{y},\boldsymbol{\xi})$, i.e., Subproblem~\eqref{SubProblem}. Let $\boldsymbol{\nu},\boldsymbol{\tau},\boldsymbol{\eta}, \boldsymbol{\psi}$ be dual variables corresponding to Constraints \eqref{SP:C1}-\eqref{SP:C4} of Subproblem~\eqref{SubProblem}. Then, the dual form of $f(\mathbf{x},\mathbf{y},\boldsymbol{\xi})$ equals:
\begin{subequations} \label{DUALSP}
\begin{align}
f(\mathbf{x},\mathbf{y},\boldsymbol{\xi})=\max_{\boldsymbol{\nu},\boldsymbol{\tau},\boldsymbol{\eta},\boldsymbol{\psi}} & \quad  \sum_{t \in \mathcal{T}} \left( \sum_{i \in \mathcal{N}^S}  \sum_{t' \in [1,t']} - y_{it'} \nu_{it} + \sum_{j \in \mathcal{N}^D} \left( \xi_{jt} \psi_{jt} - R_{jt} \xi_{jt} \right) \right), \span \span \\
\text{s.t.} & \quad - \nu_{it} + \tau_{it} \leq c^{P}_{it}, \quad &\forall&  i \in \mathcal{N}^S, t \in \mathcal{T}, \label{DUALSP:C1} \\
& \quad \eta_{it}  \leq c^{I}_{it}, \quad &\forall&  i \in \mathcal{N}^I, t \in \mathcal{T} \label{DUALSP:C2},  \\
& \quad - \tau_{it} +  \psi_{jt}  \leq c^{T}_{ijt},  \quad &\forall&  i \in \mathcal{N}^S, j \in \mathcal{N}^D, t \in \mathcal{T}, \label{DUALSP:C3}  \\
& \quad - \eta_{it}  + \psi_{jt}  \leq c^{T}_{ijt},  \quad &\forall&  i \in \mathcal{N}^I, j \in \mathcal{N}^D, t \in \mathcal{T}, \label{DUALSP:C4}  \\
& \quad \boldsymbol{\nu} \geq 0. \label{DUALSP:C5}
\end{align}\end{subequations}

Given that the subproblem has relatively complete recourse, we apply the Karush–Kuhn–Tucker (KKT) conditions as in \cite{zeng2013solving} to derive the following equivalent formulation for Model~\eqref{DUALSP}.
\begin{subequations} 
\begin{align}
\max_{\mathbf{h},\mathbf{z},\mathbf{v},\boldsymbol{\nu},\boldsymbol{\tau},\boldsymbol{\eta},\boldsymbol{\psi},\boldsymbol{\xi}} & \quad  \sum_{t \in \mathcal{T}} \left( \sum_{i \in \mathcal{N}^S} c^{P}_{it}h_{it}  + \sum_{i \in \mathcal{N}^I}  c^{I}_{it} v_{it} + \sum_{j \in \mathcal{N}^D} \left(  \sum_{i \in \mathcal{N}^S \cup \mathcal{N}^I } c^{T}_{ijt}   z_{ijt}  -  R_{jt} \xi_{jt}  \right)  \right), \span  \span \\  
\text{s.t.}&\quad \eqref{SP:C1}-\eqref{SP:C5}, \eqref{DUALSP:C1}-\eqref{DUALSP:C5}, \nonumber \\
& \quad  \Big(\sum_{t' \in [1,t]} y_{it'} - h_{it}\Big) \nu_{it} = 0,    \hspace{2cm} &\forall& i \in \mathcal{N}^S, t \in \mathcal{T}, \label{CSP1}   \\ 
& \quad (c^{P}_{it} + \nu_{it} - \tau_{it}) h_{it}= 0, \quad  &\forall&  i \in \mathcal{N}^S, t \in \mathcal{T},  \label{CSD1}  \\
& \quad ( c^{I}_{it} - \eta_{it}) v_{it}  = 0, \quad &\forall& i \in \mathcal{N}^I, t \in \mathcal{T},   \label{CSD2} \\
& \quad  (c^{T}_{ijt}  + \tau_{it} - \psi_{jt}) z_{ijt} = 0, 
 \quad &\forall&  i \in \mathcal{N}^S, j \in \mathcal{N}^D, t \in \mathcal{T},   \label{CSD3}   \\
& \quad (c^{T}_{ijt}  + \eta_{it} - \psi_{jt})  z_{ijt} = 0,   \quad &\forall&  i \in \mathcal{N}^I, j \in \mathcal{N}^D, t \in \mathcal{T}, \label{CSD4}  \\
&  \quad \xi_{jt} \in [\overline{\xi}_{jt}, \underline{\xi}_{jt}], \quad &\forall& j \in \mathcal{N}^D, t \in \mathcal{T}. \label{CSBounds}
\end{align}\end{subequations}

The above formulation includes the constraints of the Primal Subproblem~\eqref{SubProblem} and the Dual Subproblem~\eqref{DUALSP}, together with Constraints~\eqref{CSP1}--\eqref{CSBounds}, which correspond to complementary slackness conditions. We note that there is no need to include complementary slackness conditions for equality Constraints~\eqref{SP:C2}-\eqref{SP:C4} in Primal Subproblem~\eqref{SubProblem}, as their slack variables are always zero. Complementary Slackness Constraints~\eqref{CSP1}-\eqref{CSD4} include nonlinear terms. For linearization purposes, we denote $M$ as a large number and introduce a set of new binary variables $\boldsymbol{\kappa}^m$ for $m = \{1,\dots,5\}$, each corresponding to Constraints~\eqref{CSP1}-\eqref{CSD4}. Then, an equivalent set of linear constraints can be obtained as follows:
\begin{subequations}
\begin{align}
&  ( \sum_{t' \in [1,t]} y_{it'} - h_{it} ) \leq M \kappa^1_{it}, & \nu_{it} & \leq M(1-\kappa^1_{it}),   &\forall & i \in \mathcal{N}^S, t \in \mathcal{T}, \label{Linear_CS1}\\
&  (c^{P}_{it} + \nu_{it} - \tau_{it}) \leq M \kappa^2_{it}, & h_{it} & \leq M(1-\kappa^2_{it}),    & \forall& i \in \mathcal{N}^S, t \in \mathcal{T},  \\
& ( c^{I}_{it} - \eta_{it}) \leq  M \kappa^3_{it}, &v_{it}& \leq M(1-\kappa^3_{it}),   &\forall & i \in \mathcal{N}^I, t \in \mathcal{T},    \\
& (c^{T}_{ijt}  + \tau_{it} - \psi_{jt}) \leq M \kappa^4_{ijt}, &z_{ijt}& \leq M(1-\kappa^4_{ijt}), 
 &\forall&  i \in \mathcal{N}^S, j \in \mathcal{N}^D, t \in \mathcal{T},      \\
& (c^{T}_{ijt}  + \eta_{it} - \psi_{jt}) \leq M \kappa^5_{ijt}, &z_{ijt}& \leq M(1-\kappa^5_{ijt}),     &\forall & i \in \mathcal{N}^I, j \in \mathcal{N}^D, t \in \mathcal{T},  \label{Linear_CS5} \\&\kappa^1_{it},\kappa^2_{it},\kappa^3_{it},\kappa^4_{ijt},\kappa^5_{ijt}  \in \{0, 1\},   & & & \forall &  i \in \mathcal{N}, j \in \mathcal{N}^D, t \in \mathcal{T}. \label{Linear_CSbound}
\end{align} 
\end{subequations}

\begin{thm} An equivalent formulation to Model~\eqref{RobustModel} under the location-based moment function~\eqref{Decision-dependent moment} is given by:
\begin{subequations} \label{Continuous Reformulation}
\begin{align}
\min_{\substack{\mathbf{x},\mathbf{y},\alpha,\\ \boldsymbol{\beta}^1, \boldsymbol{\beta}^2, \boldsymbol{\Phi}}} & \quad  \sum_{t \in T} \sum_{i \in \mathcal{N}^S} \left( c^{F}_{it} (x_{it}- x_{i(t-1)} )+ c^{V}_{it} y_{it} \right)  + \alpha + \sum_{j \in \mathcal{N}^D} \sum_{t \in \mathcal{T}} \left(  \bar \mu_{jt} \beta^1_{jt} + \bar \mu_{jt} \sum_{i \in \mathcal{N}^S} \lambda^{l}_{ij} \Phi^1_{ijt} \right. \nonumber \\
& \hspace{6cm} \left. + \epsilon_{jt} \beta^1_{jt} -  \bar \mu_{jt} \beta^2_{jt} - \bar \mu_{jt} \sum_{i \in \mathcal{N}^S} \lambda^{l}_{ij} \Phi^2_{ijt}  + \epsilon_{jt}\beta^2_{jt} \right), \label{Continuous Reformulation:obj}\\
\text{s.t.} & \quad \eqref{RobustModel_C1}-\eqref{RobustModel_C3}, \eqref{DI: C2}, \eqref{MC1}-\eqref{MC4},  \nonumber  \\ 
& \quad \alpha \geq \max_{\boldsymbol{\xi}\in \mathcal{S}} \left\{ \sum_{t \in \mathcal{T}} \left( \sum_{i \in \mathcal{N}^S} c^{P}_{it}h_{it}  + \sum_{i \in \mathcal{N}^I}  c^{I}_{it} v_{it} + \sum_{j \in \mathcal{N}^D} \left(
  \sum_{i \in \mathcal{N}^S \cup \mathcal{N}^I} c^{T}_{ijt}   z_{ijt} -  R_{jt} \xi_{jt} \right) \right. \right. \nonumber \\ 
&\hspace{9cm}  \left. \left.  -  \sum_{j \in \mathcal{N}^D} (\beta^1_{jt} -\beta^2_{jt}) \xi_{jt} \right) \right\}, \label{Continuous Cut} \\
&\quad \eqref{SP:C1}-\eqref{SP:C5}, \eqref{DUALSP:C1}-\eqref{DUALSP:C5}, \eqref{CSBounds}, \eqref{Linear_CS1}-\eqref{Linear_CSbound}. \nonumber 
\end{align}\end{subequations}
\end{thm}

\begin{rmrk}
If the support set $\mathcal{S}$ in Model~\eqref{ReformulationFirstStep} is discrete, instead of Reformulation~\eqref{Continuous Reformulation}, one can derive a single-level mixed-integer linear programming reformulation, which can be directly solved using off-the-shelf optimization solvers.  
More specifically, let  $\mathcal{K} = \{\boldsymbol{\xi}^1,\dots, \boldsymbol{\xi}^K\}$ be a discrete support set where $\mathcal{K} \subseteq \mathcal{S}$. In this case, each $\boldsymbol{\xi}^k$ represents a matrix of dimensions $J \times T$, containing values $\xi_{jt}$ corresponding to demand locations $j \in \mathcal{N}^D$ and times $t \in \mathcal{T}$. This alternative reformulation under a support set $\mathcal{K}$ is presented in Theorem~\ref{thm:DiscreteReformulation}.

\begin{thm} \label{thm:DiscreteReformulation}
An equivalent single-stage formulation to the original Model~\eqref{RobustModel} under support set $\mathcal{K} \subseteq \mathcal{S}$ and location-based moment function~\eqref{Decision-dependent moment}  can be obtained as follows:
\begin{subequations} \label{Discrete Reformulation}
\begin{align}
\min_{\substack{\mathbf{x},\mathbf{y},\alpha,\\ \boldsymbol{\beta}^1, \boldsymbol{\beta}^2, \boldsymbol{\Phi}}} & \quad  \sum_{t \in T} \sum_{i \in \mathcal{N}^S} \left( c^{F}_{it} (x_{it}- x_{i(t-1)} )+ c^{V}_{it} y_{it} \right)  + \alpha + \sum_{j \in \mathcal{N}^D} \sum_{t \in \mathcal{T}} \left(  \bar \mu_{jt} \beta^1_{jt} + \bar \mu_{jt} \sum_{i \in \mathcal{N}^S} \lambda^{l}_{ij} \Phi^1_{ijt} \right. \span \span \span  \nonumber \\
& \hspace{6cm} \left. + \epsilon_{jt} \beta^1_{jt} -  \bar \mu_{jt} \beta^2_{jt} - \bar \mu_{jt} \sum_{i \in \mathcal{N}^S} \lambda^{l}_{ij} \Phi^2_{ijt}  + \epsilon_{jt}\beta^2_{jt} \right), \span \span \span \\
\text{s.t.} &  \quad \eqref{RobustModel_C1}-\eqref{RobustModel_C3}, \eqref{DI: C2}, \eqref{MC1}-\eqref{MC4}, \span \span  \span \nonumber \\
& \quad   \alpha + \sum_{j \in \mathcal{N}^D} \sum_{t \in \mathcal{T}}  (\beta^1_{jt} -\beta^2_{jt}) \xi^k_{jt}   \geq \sum_{t \in \mathcal{T}} \left( \sum_{i \in \mathcal{N}^S} c^{P}_{it}h_{it}^k  + \sum_{i \in \mathcal{N}^I}  c^{I}_{it} v_{it}^k  \right. \span \span \span   \nonumber \\
& \hspace{5.4cm} \left. + \sum_{j \in \mathcal{N}^D} \left( \sum_{i \in \mathcal{N}^S \cup \mathcal{N}^I} c^{T}_{ijt}   z_{ijt}^k -  R_{jt} \xi^k_{jt} \right) \right), \quad \forall k \in \mathcal{K}, \span \span \span \label{Discrete_Calpha} \\
& \quad   h_{it}^k  \leq \sum_{t' \in [1,t]} y_{it'},  \hspace{3cm} &\forall& i \in \mathcal{N}^S, t \in \mathcal{T},  k \in \mathcal{K},  \\ 
& \quad h^k_{it}  =  \sum_{j \in \mathcal{N}} z^k_{ijt}, \quad &\forall& i \in \mathcal{N}^S, t \in \mathcal{T}, k \in \mathcal{K}, \\ 
&  \quad  v^k_{it}  =  \sum_{j \in \mathcal{N}} z^k_{ijt}, \quad &\forall& i \in \mathcal{N}^I, t \in \mathcal{T}, k \in \mathcal{K}, \\ 
& \quad   \sum_{i \in \mathcal{N}} z^k_{ijt} = \xi^k_{jt},  \quad &\forall& j \in \mathcal{N}^D, t \in \mathcal{T},  k \in \mathcal{K}, \\ 
& \quad   \mathbf{h},\mathbf{z},\mathbf{v} \geq 0.  \label{Discrete_domain}
\end{align}\end{subequations}
\end{thm}
\end{rmrk}

We note that Formulation~\eqref{Discrete Reformulation} remains valid for continuous support sets; however, it transforms into a semi-infinite problem as Constraints~\eqref{Discrete_Calpha}-\eqref{Discrete_domain} are infinitely many.  Consequently, Formulation~\eqref{Discrete Reformulation} becomes computationally intractable for these sets, requiring further solution approaches. 

\section{Solution Algorithm} \label{sec:SolutionApproach} 
In this section, we introduce a tailored algorithm to solve the proposed distributionally robust decision-dependent hydrogen network expansion planning problem. Our algorithm is an enhanced version of the C\&CG algorithm initially proposed by \cite{zeng2013solving} to solve two-stage robust optimization problems. The algorithm employs a decomposition strategy, iteratively introducing columns and constraints tailored to worst-case scenarios identified during the optimization process. Through this iterative process, the algorithm avoids the computational burden of considering the entire scenario space simultaneously, enhancing the efficiency of the optimization process. In what follows, we describe the classical C\&CG algorithm, introduce our algorithm acceleration strategies, and present a row generation algorithm to be used as a benchmark in our computational experiments.

Model~\eqref{Continuous Reformulation} cannot be directly solved by commercial optimization solvers as it contains an inner maximization problem in Constraint~\eqref{Continuous Cut}. To address this computational complexity, we decompose Model~\eqref{Continuous Reformulation} into a master and subproblem to obtain lower bound (LB) and upper bound (UB) for the original problem. Initially, the master problem includes the first-stage variables together with a subset of the scenarios $\mathcal{S}' \subseteq \mathcal{S}$. This relaxed formulation acts as a starting point and offers a valid LB for the original Model~\eqref{Continuous Reformulation}. Accordingly, we define the master problem $(MP)$ as follows:
\begin{subequations} 
\begin{align}
\min_{\substack{\mathbf{x},\mathbf{y},\alpha,\\ \boldsymbol{\beta}^1, \boldsymbol{\beta}^2, \boldsymbol{\Phi}}} & \quad  \sum_{t \in T} \sum_{i \in \mathcal{N}^S} \left( c^{F}_{it} (x_{it}- x_{i(t-1)} )+ c^{V}_{it} y_{it} \right)  + \alpha + \sum_{j \in \mathcal{N}^D} \sum_{t \in \mathcal{T}} \left(  \bar \mu_{jt} \beta^1_{jt} + \bar \mu_{jt} \sum_{i \in \mathcal{N}^S} \lambda^{l}_{ij} \Phi^1_{ijt} \right. \span \span \nonumber \\
& \hspace{6cm} \left. + \epsilon_{jt} \beta^1_{jt} -  \bar \mu_{jt} \beta^2_{jt} - \bar \mu_{jt} \sum_{i \in \mathcal{N}^S} \lambda^{l}_{ij} \Phi^2_{ijt}  + \epsilon_{jt}\beta^2_{jt} \right), \span \span  \label{Master: obj} \\
\text{s.t.} & \quad \eqref{RobustModel_C1}-\eqref{RobustModel_C3}, \eqref{DI: C2}, \eqref{MC1}-\eqref{MC4}, \span \span  \nonumber  \\
& \quad \alpha \geq \sum_{t \in \mathcal{T}} \left( \sum_{i \in \mathcal{N}^S} c^{P}_{it}h_{it}  + \sum_{i \in \mathcal{N}^I}  c^{I}_{it} v_{it} + \sum_{j \in \mathcal{N}^D} \left(
  \sum_{i \in \mathcal{N}^S \cup \mathcal{N}^I} c^{T}_{ijt}   z_{ijt} -  R_{jt} \xi^k_{jt}  \right) \right. \span \span \nonumber \\
& \hspace{8cm} \left.  -  \sum_{j \in \mathcal{N}^D} (\beta^1_{jt} -\beta^2_{jt}) \xi^k_{jt}  \right), \quad k \in \mathcal{S}', \span \span  \label{Master: cut}\\
& \quad   h_{it}^k  \leq \sum_{t' \in [1,t]} y_{it'},  \hspace{3cm}  &\forall& i \in \mathcal{N}^S, t \in \mathcal{T},  k \in \mathcal{S}',  \\ 
& \quad h^k_{it}  =  \sum_{j \in \mathcal{N}} z^k_{ijt}, \quad &\forall& i \in \mathcal{N}^S, t \in \mathcal{T}, k \in \mathcal{S}', \\ 
&  \quad  v^k_{it}  =  \sum_{j \in \mathcal{N}} z^k_{ijt}, \quad &\forall& i \in \mathcal{N}^I, t \in \mathcal{T}, k \in \mathcal{S}', \\ 
& \quad   \sum_{i \in \mathcal{N}} z^k_{ijt} = \xi^k_{jt},  \quad &\forall& j \in \mathcal{N}^D, t \in \mathcal{T},  k \in \mathcal{S}', \\ 
& \quad   \mathbf{h},\mathbf{z},\mathbf{v} \geq 0.  
\end{align}\end{subequations}

The important aspect is to create a subset $\mathcal{S}'$ that includes critical scenarios $\boldsymbol{\xi}^k$ to effectively approximate the original model. To achieve this, we solve the below Subproblem $(SP)$~\eqref{SP:CCGA} under given first-stage variables $\mathbf{x},\mathbf{y},\boldsymbol{\beta}^1, \boldsymbol{\beta}^2$  to identify the worst-case scenarios $\boldsymbol{\xi}^*$ to be added to $\mathcal{S}'$.
\begin{subequations} \label{SP:CCGA}
    \begin{align}
\max_{\boldsymbol{\xi}\in \mathcal{S}} & \quad \sum_{t \in \mathcal{T}} \left( \sum_{i \in \mathcal{N}^S} c^{P}_{it}h_{it}  + \sum_{i \in \mathcal{N}^I}  c^{I}_{it} v_{it} + \sum_{j \in \mathcal{N}^D} \left(
  \sum_{i \in \mathcal{N}^S \cup \mathcal{N}^I} c^{T}_{ijt}   z_{ijt} -  R_{jt} \xi_{jt} \right) -  \sum_{j \in \mathcal{N}^D} (\beta^1_{jt} -\beta^2_{jt}) \xi_{jt} \right), \\
\text{s.t.} &\quad \eqref{SP:C1}-\eqref{SP:C5}, \eqref{DUALSP:C1}-\eqref{DUALSP:C5},\eqref{CSBounds}, \eqref{Linear_CS1}-\eqref{Linear_CSbound}. \nonumber 
\end{align}\end{subequations}

We update the UB based on the solution obtained by the master problem and the optimal objective function value of the subproblem under this solution. By gradually expanding the subset $\mathcal{S}'$ with the inclusion of the new scenarios obtained from the subproblem, we can achieve stronger LBs for the master problem. The algorithm stops when the optimality gap $\%\text{Gap} = 100\times(1-\text{LB}/\text{UB})$ reaches a predetermined threshold. 

\subsection{Algorithm Acceleration Strategies}
To improve the convergence of the algorithm, we propose three different acceleration strategies by solving the master problem to a relative optimality gap, including pre-generated columns and constraints, and parallelizing the subproblem solution procedure.

The master problem's solution time is a major bottleneck in the efficiency of the C\&CG algorithm. Through the course of the algorithm, the master problem grows substantially in size and requires to be solved multiple times. Therefore, instead of solving the master problem to optimality in each iteration, we solve inexact master problems by employing a relative optimality gap ($\delta\%$). We start with an initial relative gap for the master problem and update it based on the difference between the lower and upper bounds obtained in each iteration. We ensure that the relative optimality gap of the master problem converges to a very small threshold. We establish both an initial ($\delta^i$) and a final ($\delta^f$) relative optimality gap. Introducing a factor ($f^{\%}$), which ranges between zero and one, enables us to fine-tune the relative optimality gap as a proportion of the optimality gap (i.e.,\%\text{Gap}). The formulation of the relative optimality gap is as follows:
\begin{equation}
\delta\%=\max \left\{\min \left\{\delta^i, f^{\%} \times \frac{UB-LB}{UB}\right\}, \delta^f \right\} \label{eq:Relative gap}
\end{equation}
Equation~\eqref{eq:Relative gap} ensures that the relative optimality gap is bounded by the initial and final thresholds while being dynamically adjusted based on the difference between upper and lower bounds.

We note that Subproblem \eqref{SP:CCGA}, denoted as $(SP)$, can be decomposed into time periods when the first-stage decisions are fixed. We utilize this property and decompose the subproblems into T subproblems, denoted by $(SP_t)$, to solve them in parallel. In addition, we can generate Constraint~\eqref{Master: cut} for each subproblem. Accordingly, we define $\alpha_t$ for each $t \in \mathcal{T}$ and replace Constraint~\eqref{Master: cut} with the following constraint.
\begin{equation}
\alpha_t \geq \sum\limits_{i \in \mathcal{N}^S} c^{P}_{it}h_{it} + \sum\limits_{i \in \mathcal{N}^I} c^{I}_{it} v_{it} + \sum\limits_{j \in \mathcal{N}^D} \left(
\sum\limits_{i \in \mathcal{N}^S \cup \mathcal{N}^I} c^{T}_{ijt} z_{ijt} - R_{jt} \xi^k_{jt} \right) - \sum\limits_{j \in \mathcal{N}^D} (\beta^1_{jt} -\beta^2_{jt}) \xi^k_{jt}, \quad k \in \mathcal{S}', t \in \mathcal{T}.
\end{equation}
We then replace the $\alpha$ term in Objective Function~\eqref{Master: obj} with $\sum_{t \in \mathcal{T}} \alpha_t$. This approach makes the cuts more specific, potentially leading to improved lower bounds.

To further improve the computational performance, we generate columns and constraints for a predetermined set of scenarios at the beginning of our algorithm. In our preliminary experiments,  we observe that the worst-case scenarios identified by the classical algorithm often exhibit a specific pattern: in these scenarios, typically, only one decision variable $\xi_{jt}$ reaches its upper bound value $\overline{\xi}_{jt}$, while the others are at their lower bound values $\underline{\xi}_{jt}$. Leveraging these patterns, we generate $|J|*|T|$ distinct cases that can be potential worst cases of the problem. In each of these potential cases, only one $\xi_{jt}$ takes its upper bound value $\overline{\xi}_{jt}$, while the others take their lower bound values $\underline{\xi}_{jt}$. We then add this predetermined set of worst-case scenarios to set $\mathcal{S}'$ at the first step of the algorithm so that corresponding columns and constraints are incorporated into the master problem from the beginning.

\begin{algorithm}[h]
\caption{Solution Algorithm}
\label{SolutionAlgorithm}
\begin{algorithmic}[1]
\STATE Set LB = -$\infty$, UB =$\infty$, Iter=0, $\delta\% = \delta^i$.
\STATE Initialize $\mathcal{S}'$ with predetermined worst-case scenarios.
\WHILE {$100\times(1-\text{LB}/\text{UB})<\epsilon $}
\STATE Solve $(MP)$ to relative optimality gap $\delta\%$ to obtain $z(MP)$ and $\mathbf{x}^*,\mathbf{y}^*,\boldsymbol{\beta}^{1*}, \boldsymbol{\beta}^{2*}$.
\STATE Set LB $\leftarrow z(MP)$.
\STATE Fix $\mathbf{x}^*,\mathbf{y}^*,\boldsymbol{\beta}^{1*}, \boldsymbol{\beta}^{2*}$ and decompose $(SP)$ into $|\mathcal{T}|$ subproblems.
\FORALL{$t \in \mathcal{T}\ (\text{in parallel})$} 
\STATE Solve $(SP_t)$  to obtain  $z(SP_t)$ and $\boldsymbol{\xi_t}^*$.
\ENDFOR
\STATE  \text{UB}$^\text{c}$ = $\sum_{t \in T} \sum_{i \in \mathcal{N}^S} \left( c^{F}_{it} (x_{it}^*- x_{i(t-1)}^* )+ c^{V}_{it} y_{it}^* \right)  + \sum_{j \in \mathcal{N}^D} \sum_{t \in \mathcal{T}} \left(  \bar \mu_{jt} \beta^{1*}_{jt} + \bar \mu_{jt} \sum_{i \in \mathcal{N}^S} \lambda^{l}_{ij} \Phi^{1*}_{ijt} \right. \hspace{2cm} $  \\
$\left. \hspace{1.1cm} + \epsilon_{jt} \beta^{1*}_{jt} -  \bar \mu_{jt} \beta^{2*}_{jt} - \bar \mu_{jt} \sum_{i \in \mathcal{N}^S} \lambda^{l}_{ij} \Phi^{2*}_{ijt}  + \epsilon_{jt}\beta^{2*}_{jt} \right) + \sum_{t \in \mathcal{T}}z(SP_t)$. \label{step_1}
\IF{\text{UB}$^\text{c}$ $\leq $ UB}  
\STATE UB $\leftarrow$ \text{UB}$^\text{c}$.
\ENDIF \label{step_2}
\STATE $\mathcal{S}' \leftarrow \mathcal{S}' \cup \{\boldsymbol{\xi}^*\}$. 
\STATE Generate columns and constraints for $\boldsymbol{\xi}^*$ and add them to $(MP)$. 
\STATE $\delta\% \leftarrow \max \left\{\min \left\{\delta^i, f^{\%} \times \frac{UB-LB}{UB}\right\}, \delta^f \right\}$.
\STATE $\%\text{Gap} \leftarrow 100\times(1-\text{LB}/\text{UB})$. 
\STATE Iter $\leftarrow$ Iter + 1.
\ENDWHILE 
\end{algorithmic}
\end{algorithm}

We summarize our solution approach in Algorithm~\ref{SolutionAlgorithm}. The algorithm starts with incorporating predetermined worst-case scenarios to initialize the set $\mathcal{S}'$. Subsequently, in each iteration, the master problem $(MP)$ is solved to the given relative optimality gap $\delta^{\%}$ to derive an optimal objective value $z(MP)$ and the solution vectors $\mathbf{x}^*, \mathbf{y}^*, \boldsymbol{\beta}^{1*}, \boldsymbol{\beta}^{2*}$. The LB is then updated with the obtained objective value.  After fixing the variables $\mathbf{x}^*,\mathbf{y}^*,\boldsymbol{\beta}^{1*}, \boldsymbol{\beta}^{2*}$ in $(SP)$, the remaining problem is decomposed into a set of subproblems $(SP_t)$, which are solved in parallel.  Upon obtaining the optimal objective values $z(SP_t)$ and the corresponding worst-case scenarios $\boldsymbol{\xi_t}^*$ for each subproblem, the UB is updated if the candidate upper bound  \text{UB}$^\text{c}$ is better. In addition, the identified worst-case scenario $\boldsymbol{\xi}^*$ is incorporated into the set $\mathcal{S}'$ to generate corresponding columns and constraints. This iterative process continues until the optimality gap \%Gap diminishes below the specified threshold~$\epsilon$.

\begin{rmrk}
In the case of discrete support sets, Model~\eqref{Discrete Reformulation} can be directly solved using commercial optimization solvers. However, computational challenges arise when dealing with large support sets. For such cases our solution algorithm can be used by changing Subproblem ~\eqref{SP:CCGA} with the following subproblem, while   the rest of the solution algorithm remains unchanged. 
\end{rmrk}
\begin{align}
\max_{\boldsymbol{\xi}\in \mathcal{K}} & \quad \sum_{t \in \mathcal{T}} \left( \sum_{i \in \mathcal{N}^S} c^{P}_{it}h_{it}  + \sum_{i \in \mathcal{N}^I}  c^{I}_{it} v_{it} + \sum_{j \in \mathcal{N}^D} \left(
  \sum_{i \in \mathcal{N}^S \cup \mathcal{N}^I} c^{T}_{ijt}   z_{ijt} -  R_{jt} \xi_{jt} \right) -  \sum_{j \in \mathcal{N}^D} (\beta^1_{jt} -\beta^2_{jt}) \xi_{jt} \right),  \label{SP:Discrete} \\
\text{s.t.} &\quad \eqref{SP:C1}-\eqref{SP:C5}. \nonumber 
\end{align}

\subsection{Benders Decomposition as a Row Generation Algorithm}
Alternatively, we consider applying Benders decomposition as a benchmark algorithm in our computational experiments to evaluate the performance our proposed approach. Benders decomposition is a well-known row generation algorithm for solving various optimization problems, including robust optimization. This row generation algorithm is conceptually similar to our row-and-column generation algorithm in terms of  decomposing the problem into a master and subproblem and iteratively adding cuts obtained from the subproblem to the master problem. We define the master problem $(BMP)$ for the Benders decomposition algorithm as follows:
\begin{subequations} 
\begin{align}
\min_{\substack{\mathbf{x},\mathbf{y},\alpha,\\ \boldsymbol{\beta}^1, \boldsymbol{\beta}^2, \boldsymbol{\Phi}}} & \quad  \sum_{t \in T} \sum_{i \in \mathcal{N}^S} \left( c^{F}_{it} (x_{it}- x_{i(t-1)} )+ c^{V}_{it} y_{it} \right)  + \alpha + \sum_{j \in \mathcal{N}^D} \sum_{t \in \mathcal{T}} \left(  \bar \mu_{jt} \beta^1_{jt} + \bar \mu_{jt} \sum_{i \in \mathcal{N}^S} \lambda^{l}_{ij} \Phi^1_{ijt} \right. \span \span \nonumber \\
& \hspace{6cm} \left. + \epsilon_{jt} \beta^1_{jt} -  \bar \mu_{jt} \beta^2_{jt} - \bar \mu_{jt} \sum_{i \in \mathcal{N}^S} \lambda^{l}_{ij} \Phi^2_{ijt}  + \epsilon_{jt}\beta^2_{jt} \right), \span \span   \\
\text{s.t.} & \quad \eqref{RobustModel_C1}-\eqref{RobustModel_C3}, \eqref{DI: C2}, \eqref{MC1}-\eqref{MC4}, \span \span  \nonumber  \\
& \quad   \alpha  \geq  \sum_{t \in \mathcal{T}} \left( \sum_{i \in \mathcal{N}^S}  \sum_{t' \in [1,t']} - y_{it'} \nu_{it}^* + \sum_{j \in \mathcal{N}^D}  \left( \xi_{jt}^*  \psi_{jt}^* -  R_{jt}\xi_{jt}^*  - (\beta^1_{jt} -\beta^2_{jt}) \xi_{jt}^* \right)  \right) \quad \boldsymbol{\xi}^* \in \Theta \label{BendersCut}
\end{align}\end{subequations}

Here, constraint \eqref{BendersCut} corresponds to the Benders cut which is obtained based on the dual of the subproblem given in Model \eqref{DUALSP}. At the beginning of the algorithm, we start with an empty set of $\Theta$.  Then, in each iteration, we solve $(SP)$ (model \eqref{SubProblem}) to obtain a worst-case scenario $\boldsymbol{\xi}^*$ and optimal values of the dual variables $\boldsymbol{\nu}^*,\boldsymbol{\psi}^*$ to generate the cut in the form of constraint \eqref{BendersCut}.  For each iteration, the lower bound (LB) is obtained as the objective value of the Benders master problem $z(BMP)$. The upper bound (UB) is computed according to steps \ref{step_1} and \ref{step_2} of Algorithm \ref{SolutionAlgorithm}. The algorithm iterates until $\%\text{Gap}$ reaches the predetermined threshold $\epsilon$.

\section{Computational Analyses}\label{sec:ComputationalExperiments}
This section provides a comprehensive set of results highlighting the computational efficiency of the provided solution approach, and the findings and practical implications of the proposed optimization model under a real case study. 
In Section \ref{sec:Computational}, 
we provide a computational study demonstrating performance of our model reformulations and the enhanced C\&CG algorithm, where we structurally compare and evaluate our model and algorithm components for several benchmark instances for the hydrogen generation expansion planning problem. Afterwards, in Section \ref{sec:CaseStudy}, we present a real case study illustrating the value of our model and associated 
optimization frameworks 
to support a cost-efficient transition towards a green hydrogen economy in the Netherlands as part of the project HEAVENN - Hydrogen Energy Applications in Valley Environments for Northern Netherlands \citep{Heavenn}.

\subsection{Computational Study}
\label{sec:Computational}
For  illustrating the computational performance of our solution approaches, we generate an extensive set of random instances with different numbers of supply and demand nodes. We construct sets of supply $\mathcal{N}^S$, import $\mathcal{N}^I$, and demand $\mathcal{N}^D$ nodes within a 100x100 coordinate system. The set-up cost $c^{F}_{i1}$ and variable cost $c^{VS}_{i1}$ of building an electrolyzer at the first planning period to location $i \in \mathcal{N}^S$ are derived randomly from normal distributions $\mathcal{N}(4000,400)$ and $\mathcal{N}(20,2)$, respectively. We assume a 50\% periodic decrease in the set-up and variable costs due to expected advancements in technology.
We set production $c^{VP}_{i1}$ and import $c^{I}_{i1}$ costs for each node $i$ at the first period to 10 and 150, respectively. We assume that the transportation cost per unit distance is 0.5. We let $\Delta_{ij}$ be the Euclidean distance between nodes $i,j \in \mathcal{N}$. Then, we set transportation cost $c^{T}_{ijt}$ between nodes $i,j$ to $\Delta_{ij}\times 0.5 $. Revenue $R_{jt}$ for supplying per demand unit is set to 100. We assume a 10\% periodic decrease in all operational costs throughout the planning period.

We sample the mean demand $\overline{\mu}_{j1}$ at each demand location $j$ at the first period from a normal distribution $\mathcal{N}(50,5)$. We assume a 50\% periodic increase  in this baseline mean demand.  We establish lower ($\underline{\xi}_{jt}$) and upper ($\overline{\xi}_{jt}$) bounds for the support set of demand at $0.75 \times \overline{\mu}_{jt}$ and $1.25 \times \overline{\mu}_{jt}$, respectively. For the computational study, we consider the location-based moment function under Equation \eqref{Decision-dependent moment}. We set the $\lambda^l_{ij}$ parameters by considering the distance $\Delta_{ij}$ between supply $i \in \mathcal{N}^S$ and demand $j \in \mathcal{N}^D$ nodes. We consider them as a decreasing function of distance $\Delta_{ij}$ to ensure that the effect of opening a facility at supply node $i$ on the demand of node $j$ is higher if they are closer to each other as in \cite{basciftci2021distributionally}. Specifically, we use the function $\exp\left(-\Delta_{ij}/25 \right)$. We assume that the opening of new supply facilities has the potential to increase baseline demand even more by up to 50\%. To accommodate this assumption, we normalize the calculated $\lambda^l_{ij}$ values for each demand location $j$ in a way that their summation ($\sum_{i \in \mathcal{N}^S} \lambda^l_{ij}$) equals 50\%.

We assess the performance of the proposed optimization frameworks by conducting computational experiments with varying model parameters. We consider instance sizes with supply nodes ranging from four to twenty, and periods ranging from one to five, including all intermediate values. For each instance size, we randomly generate three instances, leading to a total of 255 distinct instances. Other parameters remain fixed throughout the computational study unless otherwise stated. We conduct computational experiments on an Intel Xeon E5 2680v3 CPU with a 2.5 GHz processor and 32 GB RAM. The implementation is carried out in Python, utilizing \texttt{Gurobi} 10.0.0  to solve both the master problem and subproblem within the solution algorithms. A time limit of 5 hours is imposed to maintain consistency across experiments. We set the optimality gap (\%Gap) for all algorithms to 0.1\%.

\subsubsection{Computational Performance.}  
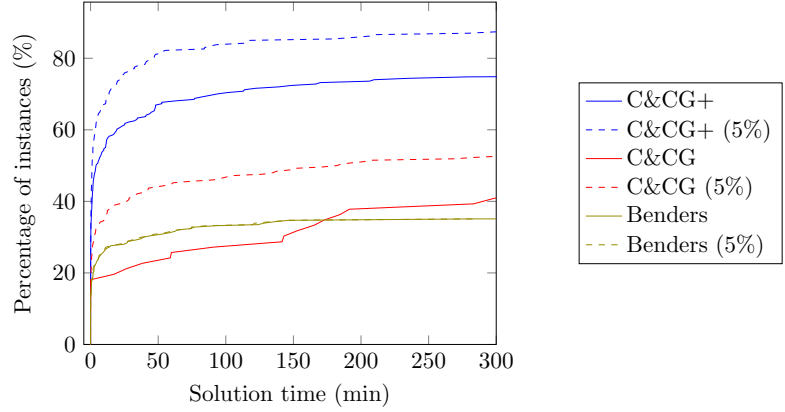
\begin{figure}[h]
\centering
\begin{tikzpicture}[scale=0.8, transform shape]
\begin{axis}[
    xlabel={Solution time (min)},
    ylabel={Percentage of instances (\%)},
    xmin=-5, xmax=300,
    ymin=0, 
    legend style={at={(1.2,0.5)}, anchor=west, legend cell align=left},
    ymajorgrids=false,
    grid style=dashed,
]
\addplot[
    color=blue,
    ]
    coordinates {
    (426.0331342, 75.10) (279.3617788, 74.80) (235.4397149, 74.40) (209.6002066, 74.00) (206.2925364, 73.60) (169.8216016, 73.20) (166.5958207, 72.80) (148.437669, 72.40) (138.9991236, 72.00) (121.8914631, 71.60) (113.7377789, 71.20) (112.9372948, 70.80) (100.8576538, 70.40) (93.36165356, 70.00) (87.87816634, 69.60) (82.35509089, 69.20) (77.08420778, 68.80) (76.53098981, 68.50) (63.21377857, 68.10) (52.72721075, 67.70) (52.62684364, 67.30) (48.07980009, 66.90) (48.00684324, 66.50) (47.80787251, 66.10) (47.35128431, 65.70) (47.03364466, 65.30) (44.84312213, 64.90) (43.83229063, 64.50) (41.10605467, 64.10) (40.69823708, 63.70) (34.70424764, 63.30) (34.21003803, 62.90) (32.32632885, 62.50) (28.4381092, 62.20) (25.4148756, 61.80) (25.11471868, 61.40) (23.13368822, 61.00) (21.55068839, 60.60) (19.54412345, 60.20) (19.21901467, 59.80) (18.33065697, 59.40) (17.69272299, 59.00) (16.56967702, 58.60) (13.73044685, 58.20) (12.79105368, 57.80) (12.14937573, 57.40) (11.78943837, 57.00) (11.66149885, 56.60) (11.57575373, 56.20) (11.31209128, 55.90) (11.28405717, 55.50) (11.07951925, 55.10) (9.786784852, 54.70) (9.743828889, 54.30) (8.699723833, 53.90) (8.112563429, 53.50) (7.819593653, 53.10) (7.479683069, 52.70) (6.768323196, 52.30) (6.474863521, 51.90) (6.062251641, 51.50) (5.572593263, 51.10) (5.405121592, 50.70) (4.307684244, 50.30) (4.032306253, 50.00) (3.864281342, 49.60) (3.717063479, 49.20) (3.534485694, 48.80) (3.441635332, 48.40) (3.091181395, 48.00) (2.932666926, 47.60) (2.627477874, 47.20) (2.543720192, 46.80) (2.367063883, 46.40) (2.099726092, 46.00) (1.988320949, 45.60) (1.975143346, 45.20) (1.83210776, 44.80) (1.772288002, 44.40) (1.73582279, 44.00) (1.706614645, 43.70) (1.699952239, 43.30) (1.550365069, 42.90) (1.429101231, 42.50) (1.405165756, 42.10) (1.287277715, 41.70) (1.282981429, 41.30) (1.16790653, 40.90) (1.140971655, 40.50) (1.096518953, 40.10) (1.059618748, 39.70) (1.058923059, 39.30) (0.939870215, 38.90) (0.910418285, 38.50) (0.90751712, 38.10) (0.86418117, 37.70) (0.726915007, 37.40) (0.655534877, 37.00) (0.639587892, 36.60) (0.568993012, 36.20) (0.565400976, 35.80) (0.537926257, 35.40) (0.499482768, 35.00) (0.450222911, 34.60) (0.380206822, 34.20) (0.322111405, 33.80) (0.297950547, 33.40) (0.289741901, 33.00) (0.251584393, 32.60) (0.225860307, 32.20) (0.204140852, 31.80) (0.203641718, 31.40) (0.188223342, 31.10) (0.182635093, 30.70) (0.179504367, 30.30) (0.176497102, 29.90) (0.174150646, 29.50) (0.165851287, 29.10) (0.159370667, 28.70) (0.148846329, 28.30) (0.141162261, 27.90) (0.136036306, 27.50) (0.13188552, 27.10) (0.121488981, 26.70) (0.120541029, 26.30) (0.119852122, 25.90) (0.112976413, 25.50) (0.107446377, 25.10) (0.098526536, 24.80) (0.097227724, 24.40) (0.095595415, 24.00) (0.091084626, 23.60) (0.080791931, 23.20) (0.079707344, 22.80) (0.07184961, 22.40) (0.070695722, 22.00) (0.069512546, 21.60) (0.063668833, 21.20) (0.058834051, 20.80) (0.058821688, 20.40) (0.050074056, 20.00) (0.049602183, 19.60) (0.049582825, 19.20) (0.041531974, 18.80) (0.039880081, 18.50) (0.038340499, 18.10) (0.037032653, 17.70) (0.030047625, 17.30) (0.030023582, 16.90) (0.029525367, 16.50) (0.026894817, 16.10) (0.024956233, 15.70) (0.022192629, 15.30) (0.020895383, 14.90) (0.01697787, 14.50) (0.014141352, 14.10) (0.012628735, 13.70) (0.012231569, 13.30) (0.011906644, 12.90) (0.010699742, 12.50) (0.010651555, 12.20) (0.010200195, 11.80) (0.009747442, 11.40) (0.006399075, 11.00) (0.006151102, 10.60) (0.006037345, 10.20) (0.005688281, 9.80) (0.005380764, 9.40) (0.00373905, 9.00) (0.002808148, 8.60) (0.00227153, 8.20) (0.002155037, 7.80) (0.001581482, 7.40) (0.001502153, 7.00) (0.001462573, 6.60) (0.00144685, 6.20) (0.001359903, 5.90) (0.00125887, 5.50) (0.00120976, 5.10) (0.001100116, 4.70) (0.001022449, 4.30) (0.00099054, 3.90) (0.000988693, 3.50) (0.000900952, 3.10) (0.000880505, 2.70) (0.000856273, 2.30) (0.000811243, 1.90) (0.000749671, 1.50) (0.000744684, 1.10) (0.000715968, 0.70) (0.000581479, 0.30) (0.000500866, 0.00)
    };

\addplot[
    color=blue,
    dashed,
    ]
    coordinates {
 (300.1151043, 87.40) (279.3617786, 87.00) (209.6002062, 86.60) (206.2925362, 86.20) (193.8274584, 85.80) (179.2174343, 85.40) (117.3238747, 85.00) (116.8150089, 84.60) (113.7377786, 84.20) (93.36165339, 83.80) (87.87816617, 83.40) (84.41902206, 83.00) (83.9853725, 82.60) (55.20055906, 82.20) (52.72721047, 81.80) (52.62684333, 81.40) (48.07979981, 81.10) (47.3512841, 80.70) (47.03364442, 80.30) (44.84312192, 79.90) (43.40703679, 79.50) (41.28939386, 79.10) (40.69823685, 78.70) (40.12766619, 78.30) (34.70424736, 77.90) (32.36752268, 77.50) (32.3263286, 77.10) (28.43810902, 76.70) (27.86331549, 76.30) (25.11471841, 75.90) (24.36656016, 75.50) (23.13368801, 75.10) (21.75203936, 74.80) (21.73723776, 74.40) (21.55068823, 74.00) (19.78131329, 73.60) (19.54412331, 73.20) (19.21901451, 72.80) (18.33065673, 72.40) (16.56967685, 72.00) (15.40810125, 71.60) (14.80571587, 71.20) (14.1751379, 70.80) (13.73044663, 70.40) (13.59299208, 70.00) (13.2269797, 69.60) (12.79105342, 69.20) (12.14937557, 68.80) (11.57575347, 68.50) (11.31209092, 68.10) (11.28405701, 67.70) (11.07951904, 67.30) (8.699723636, 66.90) (8.112563208, 66.50) (7.819593471, 66.10) (7.479682902, 65.70) (7.213529493, 65.30) (6.768322901, 64.90) (5.724562612, 64.50) (5.573783284, 64.10) (5.572593119, 63.70) (5.405121341, 63.30) (5.236078712, 62.90) (4.467878466, 62.50) (4.307684087, 62.20) (4.032305988, 61.80) (4.00841803, 61.40) (3.864281158, 61.00) (3.717063322, 60.60) (3.534485505, 60.20) (3.391432301, 59.80) (3.292163297, 59.40) (3.11011082, 59.00) (2.932666765, 58.60) (2.795413055, 58.20) (2.627477675, 57.80) (2.543719981, 57.40) (2.367063706, 57.00) (2.332455577, 56.60) (2.099725933, 56.20) (1.975143149, 55.90) (1.772287644, 55.50) (1.735822537, 55.10) (1.706614348, 54.70) (1.702940544, 54.30) (1.699951942, 53.90) (1.688702475, 53.50) (1.550364792, 53.10) (1.429101072, 52.70) (1.405165586, 52.30) (1.287277417, 51.90) (1.248415595, 51.50) (1.167906325, 51.10) (1.140971302, 50.70) (1.134274719, 50.30) (1.132148704, 50.00) (1.098770742, 49.60) (1.096518798, 49.20) (1.059618316, 48.80) (1.058922828, 48.40) (0.947055581, 48.00) (0.939869987, 47.60) (0.910418091, 47.20) (0.907516943, 46.80) (0.895914387, 46.40) (0.864180943, 46.00) (0.832313923, 45.60) (0.728571345, 45.20) (0.72691484, 44.80) (0.712655961, 44.40) (0.655534379, 44.00) (0.639587658, 43.70) (0.574350545, 43.30) (0.568992857, 42.90) (0.565400522, 42.50) (0.549071362, 42.10) (0.537926098, 41.70) (0.529445439, 41.30) (0.49948261, 40.90) (0.492230456, 40.50) (0.473459639, 40.10) (0.458480924, 39.70) (0.450222744, 39.30) (0.435367687, 38.90) (0.43103562, 38.50) (0.380206673, 38.10) (0.322111078, 37.70) (0.29795039, 37.40) (0.289741739, 37.00) (0.289508527, 36.60) (0.251584226, 36.20) (0.225860125, 35.80) (0.217788468, 35.40) (0.214910126, 35.00) (0.204140702, 34.60) (0.203641555, 34.20) (0.194392472, 33.80) (0.188223175, 33.40) (0.182634922, 33.00) (0.176496951, 32.60) (0.17457044, 32.20) (0.174150481, 31.80) (0.168776937, 31.40) (0.165851136, 31.10) (0.165588797, 30.70) (0.159370464, 30.30) (0.148846177, 29.90) (0.141162087, 29.50) (0.136036137, 29.10) (0.131885263, 28.70) (0.121488769, 28.30) (0.120797453, 27.90) (0.120540845, 27.50) (0.119851766, 27.10) (0.112976254, 26.70) (0.107446083, 26.30) (0.098526368, 25.90) (0.097227556, 25.50) (0.095595212, 25.10) (0.091084457, 24.80) (0.090331522, 24.40) (0.088970267, 24.00) (0.080791662, 23.60) (0.079707135, 23.20) (0.073553179, 22.80) (0.071849402, 22.40) (0.070695213, 22.00) (0.069512314, 21.60) (0.063668585, 21.20) (0.058833847, 20.80) (0.058821488, 20.40) (0.050073683, 20.00) (0.049601934, 19.60) (0.049582648, 19.20) (0.041531633, 18.80) (0.039879876, 18.50) (0.038340036, 18.10) (0.037032266, 17.70) (0.030047293, 17.30) (0.030023392, 16.90) (0.026894643, 16.50) (0.024956004, 16.10) (0.022192207, 15.70) (0.02089517, 15.30) (0.016977531, 14.90) (0.014141101, 14.50) (0.0126283, 14.10) (0.012230913, 13.70) (0.011906434, 13.30) (0.010699506, 12.90) (0.01065132, 12.50) (0.01062945, 12.20) (0.01019997, 11.80) (0.009747238, 11.40) (0.006398882, 11.00) (0.006150917, 10.60) (0.006037142, 10.20) (0.005688081, 9.80) (0.005380519, 9.40) (0.003738649, 9.00) (0.002807844, 8.60) (0.002271332, 8.20) (0.002154822, 7.80) (0.001580538, 7.40) (0.001501932, 7.00) (0.001461633, 6.60) (0.001446585, 6.20) (0.001359689, 5.90) (0.001258448, 5.50) (0.001209449, 5.10) (0.001099896, 4.70) (0.001022109, 4.30) (0.000990347, 3.90) (0.000988391, 3.50) (0.000900765, 3.10) (0.000880253, 2.70) (0.000855972, 2.30) (0.000810913, 1.90) (0.000749455, 1.50) (0.000744458, 1.10) (0.00071576, 0.70) (0.000581262, 0.30) (0.000500685, 0.00)
    };
\addplot[
    color=red,
    ]
    coordinates {
    (300, 40.90) (298.9191303, 40.90) (283.1871514, 39.30) (191.3263817, 37.80) (184.4581927, 36.30) (172.8986214, 34.80) (165.1864483, 33.30) (152.9456891, 31.80) (142.990341, 30.30) (141.5972228, 28.70) (90.73062828, 27.20) (59.75926959, 25.70) (58.97171215, 24.20) (38.43781617, 22.70) (26.67900092, 21.20) (17.15934015, 19.60) (0.113748345, 18.10) (0.099255995, 16.60) (0.071892374, 15.10) (0.067307725, 13.60) (0.036773122, 12.10) (0.028672483, 10.60) (0.015404859, 9.00) (0.01398686, 7.50) (0.013857119, 6.00) (0.005750223, 4.50) (0.000713291, 3.00) (0.000652393, 1.50) (0.000336388, 0.00)
    };    

    \addplot[
    color=red,
    dashed]
    coordinates {
    (306.5140698, 52.70) (284.3791559, 52.30) (269.0113981, 51.90) (208.3240996, 51.50) (202.2465639, 51.10) (191.3263813, 50.70) (184.4581923, 50.30) (182.12702, 50.00) (171.0643163, 49.60) (149.4950211, 49.20) (147.6900962, 48.80) (137.4051788, 48.40) (134.8026662, 48.00) (127.4968238, 47.60) (105.001499, 47.20) (100.3649426, 46.80) (95.82030283, 46.40) (90.73062795, 46.00) (74.18783099, 45.60) (60.90601556, 45.20) (59.75926925, 44.80) (54.5126701, 44.40) (48.05452121, 44.00) (42.88080687, 43.70) (42.52184155, 43.30) (41.114924, 42.90) (39.0107213, 42.50) (34.34961816, 42.10) (33.02447829, 41.70) (30.7527913, 41.30) (28.4356868, 40.90) (28.39118871, 40.50) (27.07976476, 40.10) (24.98049961, 39.70) (21.23019113, 39.30) (17.15933976, 38.90) (16.48125705, 38.50) (14.42535012, 38.10) (13.45134431, 37.70) (12.40451063, 37.40) (12.38553792, 37.00) (12.21668777, 36.60) (11.39599792, 36.20) (10.83869253, 35.80) (10.35947367, 35.40) (9.864431124, 35.00) (9.747459132, 34.60) (7.091700137, 34.20) (6.591302279, 33.80) (5.886874663, 33.40) (4.673215095, 33.00) (4.365817586, 32.60) (3.897666548, 32.20) (3.763062824, 31.80) (3.69677017, 31.40) (3.471623083, 31.10) (3.088570162, 30.70) (2.790947934, 30.30) (2.686455139, 29.90) (2.439806671, 29.50) (2.235020512, 29.10) (2.233444312, 28.70) (1.816915679, 28.30) (1.662513608, 27.90) (1.478389399, 27.50) (1.252185245, 27.10) (1.171073936, 26.70) (1.114996844, 26.30) (1.108774429, 25.90) (1.079550143, 25.50) (1.046060208, 25.10) (0.875354612, 24.80) (0.797931537, 24.40) (0.758319747, 24.00) (0.727152454, 23.60) (0.701038091, 23.20) (0.607554627, 22.80) (0.56294477, 22.40) (0.497025669, 22.00) (0.482732357, 21.60) (0.441183955, 21.20) (0.404648037, 20.80) (0.38140063, 20.40) (0.362922806, 20.00) (0.285445408, 19.60) (0.283351215, 19.20) (0.214615508, 18.80) (0.195040048, 18.50) (0.190442131, 18.10) (0.173985401, 17.70) (0.144958178, 17.30) (0.141055948, 16.90) (0.133897565, 16.50) (0.133113376, 16.10) (0.113748087, 15.70) (0.112012658, 15.30) (0.099255758, 14.90) (0.094033895, 14.50) (0.072526512, 14.10) (0.071892116, 13.70) (0.067827163, 13.30) (0.067307296, 12.90) (0.049057995, 12.50) (0.036772916, 12.20) (0.028672278, 11.80) (0.023754839, 11.40) (0.023298406, 11.00) (0.020201756, 10.60) (0.018702842, 10.20) (0.016817588, 9.80) (0.015404622, 9.40) (0.014211992, 9.00) (0.013986684, 8.60) (0.013856924, 8.20) (0.012787159, 7.80) (0.009616701, 7.40) (0.007011974, 7.00) (0.005750024, 6.60) (0.005719748, 6.20) (0.003427752, 5.90) (0.001115661, 5.50) (0.001082626, 5.10) (0.000714531, 4.70) (0.000713021, 4.30) (0.000706233, 3.90) (0.00065219, 3.50) (0.000651322, 3.10) (0.000635613, 2.70) (0.000567664, 2.30) (0.000502877, 1.90) (0.000336175, 1.50) (0.000240894, 1.10) (0.000211697, 0.70) (0.000174783, 0.30) (0.000166141, 0.00)
    }; 

    \addplot[
    color=olive,
    ]
    coordinates {
   (300, 35.10)  (283.4976934, 35.10) (146.7632043, 34.70) (135.9417756, 34.30) (128.1332487, 33.90) (123.9010467, 33.50) (93.06041195, 33.20) (76.50913661, 32.80) (74.38446586, 32.40) (64.08967682, 32.00) (62.24809493, 31.60) (57.85868041, 31.20) (52.05011193, 30.80) (44.57345549, 30.40) (38.95871943, 30.00) (37.75821948, 29.60) (33.48547092, 29.20) (27.56834167, 28.80) (26.98536595, 28.40) (24.62331892, 28.00) (14.99611209, 27.60) (13.78363918, 27.20) (11.52753438, 26.80) (10.83714544, 26.40) (9.890486557, 26.00) (8.990191788, 25.60) (7.836426469, 25.20) (7.190498754, 24.90) (6.341083335, 24.50) (5.795589891, 24.10) (5.550947542, 23.70) (4.671369049, 23.30) (4.24031815, 22.90) (3.812426945, 22.50) (3.356908865, 22.10) (2.561767465, 21.70) (2.558520887, 21.30) (2.517579369, 20.90) (2.509081744, 20.50) (1.907532611, 20.10) (1.894738436, 19.70) (1.738264409, 19.30) (1.692764713, 18.90) (1.501538297, 18.50) (1.062479351, 18.10) (0.771893724, 17.70) (0.687202125, 17.30) (0.644160441, 16.90) (0.623590312, 16.60) (0.593122094, 16.20) (0.524622122, 15.80) (0.477184872, 15.40) (0.413870266, 15.00) (0.351812815, 14.60) (0.318413505, 14.20) (0.080998534, 13.80) (0.05859704, 13.40) (0.058374829, 13.00) (0.052513356, 12.60) (0.043591308, 12.20) (0.03972471, 11.80) (0.032869874, 11.40) (0.023854288, 11.00) (0.01962108, 10.60) (0.01927925, 10.20) (0.018470762, 9.80) (0.017083013, 9.40) (0.014495145, 9.00) (0.013873752, 8.60) (0.009240194, 8.30) (0.009236677, 7.90) (0.008342237, 7.50) (0.005611979, 7.10) (0.005053239, 6.70) (0.003533384, 6.30) (0.002442232, 5.90) (0.000749207, 5.50) (0.000718828, 5.10) (0.000696616, 4.70) (0.000649553, 4.30) (0.000627483, 3.90) (0.000543372, 3.50) (0.000517805, 3.10) (0.000511739, 2.70) (0.000471449, 2.30) (0.000408317, 1.90) (0.000408216, 1.50) (0.000391161, 1.10) (0.000195967, 0.70) (0.000181847, 0.30) (0.00014733, 0.00)
    };  

    \addplot[
    color=olive,
    dashed
    ]
    coordinates {
   (282.4413157, 35.10) (146.7632041, 34.70) (132.5536314, 34.30) (121.3267155, 33.90) (119.8860674, 33.50) (85.96781996, 33.20) (74.54864215, 32.80) (74.38446562, 32.40) (64.0896765, 32.00) (59.29732753, 31.60) (54.93669823, 31.20) (46.91677611, 30.80) (44.57345531, 30.40) (37.7582193, 30.00) (32.07631038, 29.60) (30.90079129, 29.20) (26.8346426, 28.80) (26.22244903, 28.40) (20.4208359, 28.00) (14.99611184, 27.60) (11.39589306, 27.20) (10.59623371, 26.80) (10.06258687, 26.40) (9.890486313, 26.00) (8.601324857, 25.60) (7.108215812, 25.20) (6.341082703, 24.90) (6.233970995, 24.50) (5.795589536, 24.10) (5.550947227, 23.70) (4.671368811, 23.30) (4.126842164, 22.90) (3.8124267, 22.50) (2.509081448, 22.10) (2.477344628, 21.70) (1.939014076, 21.30) (1.907532356, 20.90) (1.894738204, 20.50) (1.738264138, 20.10) (1.692764374, 19.70) (1.673888276, 19.30) (1.501537947, 18.90) (1.119169397, 18.50) (1.06247903, 18.10) (0.771893466, 17.70) (0.660234842, 17.30) (0.644159987, 16.90) (0.593121856, 16.60) (0.524621639, 16.20) (0.524487489, 15.80) (0.41387005, 15.40) (0.356414855, 15.00) (0.338043891, 14.60) (0.318413192, 14.20) (0.080998253, 13.80) (0.058596738, 13.40) (0.058374652, 13.00) (0.052513175, 12.60) (0.043591053, 12.20) (0.03972439, 11.80) (0.032869696, 11.40) (0.023854117, 11.00) (0.019620793, 10.60) (0.019279067, 10.20) (0.018470539, 9.80) (0.017082809, 9.40) (0.014494934, 9.00) (0.013873588, 8.60) (0.009240013, 8.30) (0.009236519, 7.90) (0.008342053, 7.50) (0.005611837, 7.10) (0.005053072, 6.70) (0.003533228, 6.30) (0.002442001, 5.90) (0.000748965, 5.50) (0.000718663, 5.10) (0.000696452, 4.70) (0.000649386, 4.30) (0.000627239, 3.90) (0.000543178, 3.50) (0.000517555, 3.10) (0.000511579, 2.70) (0.000471272, 2.30) (0.000408156, 1.90) (0.000408052, 1.50) (0.000390978, 1.10) (0.000195702, 0.70) (0.000181064, 0.30) (0.00014716, 0.00)
    };  
    
    \addlegendentry{C\&CG+} 
    \addlegendentry{C\&CG+ (5\%)} 
    \addlegendentry{C\&CG} 
    \addlegendentry{C\&CG (5\%)} 
    \addlegendentry{Benders}
    \addlegendentry{Benders (5\%) } 
    
\end{axis}
\end{tikzpicture} 
\captionsetup{font=small}
\caption{Performance of C\&CG+ compared to existing approaches of the literature (C\&CG, Benders) for a 0.1\% optimality gap and a 5\% optimality gap}\label{fig:Computational Performance}
\end{figure}

We evaluate the performance of three approaches: Benders decomposition, the classical C\&CG approach, and our enhanced C\&CG+ algorithm as outlined in Algorithm~\ref{SolutionAlgorithm}. Figure~\ref{fig:Computational Performance} illustrates the computational performance in terms of solution time. Each line represents the percentage of instances solved to optimality within the given solution time. Dashed lines correspond to cases where \%Gap is set to 5\%. The results show the superiority of the C\&CG+ approach over existing approaches in the literature. C\&CG+ achieves optimality for 75\% of the instances within the time limit, with 67\% of them solved within fifty minutes. In addition, 87\% of instances are solved to 5\% optimality. Conversely, Benders Decomposition exhibits the worst performance, achieving optimality for only 35\% of instances. Even with a 5\% optimality gap, the performance of the Benders Decomposition algorithm remains unchanged, indicating its poor convergence. The classical C\&CG approach demonstrates moderate performance, solving 41\% of instances optimally and achieving a 5\% optimality gap for 53\% of instances. 

To analyze the effect of the number of supply and demand nodes (S/D) and periods $|\mathcal{T}|$ on the computational performance, we report the average results for a specific set of instances in Table~\ref{tab:NumberNodes}. The table shows the optimality gap attained within the time limit (Gap (\%)) , number of iterations (Iter.), and solution time for Benders Decomposition, C\&CG and C\&CG+. We calculate the average performance metrics based on three randomly generated instances for each given number of nodes and periods. Instances encountering memory errors during the solution of the master problem are excluded from the average calculation and marked with an asterisk. 

\begin{table}[ht]
\small
    \centering
     \caption{Computational performance under varying number of nodes and periods}\label{tab:NumberNodes}
    \begin{tabular}{p{1cm}crrrrrrrrr}
    \toprule
          \multirow{3}{1cm}{Nodes (S/D)} &  &  \multicolumn{3}{c}{Benders Decomposition} & \multicolumn{3}{c}{C\&CG  } & \multicolumn{3}{c}{C\&CG+ } \\ \cmidrule(lr){3-5} \cmidrule(lr){6-8} \cmidrule(lr){9-11}
           &  $|\mathcal{T}|$  & Gap (\%)  & Iter. & Time & Gap (\%)  & Iter. & Time &  Gap (\%) & Iter. & Time \\ \midrule
         (4/8) & 3 & $<$0.1 & 145.3 & 1.3  & $<$0.1  & 29.0 & 0.4  & $<$0.1  & 2.3  & 0.1  \\
            & 4 & $<$0.1  & 278.7 & 3.4  & $<$0.1   & 40.3 & 0.9 & $<$0.1 & 2.3  & 0.1  \\
            & 5 & $<$0.1 & 410.3 & 6.4  & $<$0.1  & 53.3 & 1.5 & $<$0.1   & 2.3   & 0.1  \\
          (6/12) & 3 & $<$0.1 & 1218.3 & 67.1  & $<$0.1  & 48.0  & 6.1  & $<$0.1 & 6.3 & 1.1  \\
            & 4 & 7.3  & 1998.0 & 191.7  & $<$0.1 & 66.0 & 13.9 & $<$0.1   &  9.0  & 3.4  \\
            & 5 & 9.4  & 2275.6 & 294.6 & $<$0.1   & 85.3 & 29.6  & $<$0.1  & 2.0   & 0.3  \\
          (8/16) & 3 & 35.6  & 1952.0 & 300.0  & $<$0.1  & 70.0  & 33.4  & $<$0.1   & 4.9  &  81.1 \\
            & 4 & 30.0  & 1852.7 & 300.0  &  $<$0.1 & 101.3 & 104.2 & $<$0.1  & 50.0   & 146.8  \\
            & 5 & 26.6  & 1838.3 & 300.0  & 0.1  & 124.0 & 238.4 & 0.1  & 46.6   & 202.3 \\
          (10/20) & 3 & 49.2  & 1350.0 & 300.0  & 0.2  & 96.0  & 250.4 & $<$0.1  & 6.7  &  16.6 \\
            & 4 & 53.3 & 1489.7  & 300.0  & 7.2 & 104.3  & 300.0  & $<$0.1  &  2.0 & 5.7  \\
            & 5 & 32.1 & 1459.6  & 300.0  & 13.1 & 100.3 & 300.0 & $<$0.1  & 8.7  & 58.5  \\
          (12/24) & 3 & 49.4  & 1103.0 & 300.0  & 17.7  & 75.0  & 300.0  & 0.2  & 8.3  & 111.5  \\
            & 4 & 53.3 & 1015.0 & 300.0  & 25.3  & 72.0 & 300.0 &  $<$0.1 & 5.6  & 87.3  \\
            & 5 & 41.2  & 1044.3 & 300.0  & 35.1  & 75.0 & 300.0 & 0.2 & 6.7  & 114.8  \\
          (14/28) & 3 & 50.4  & 684.3 & 300.0   & 48.3  & 52.6 & 300.0 & 1.7  & 6.6  & 109.8  \\
            & 4 & 53.1  & 459.7 & 300.0  & 49.1  & 54.0 & 300.0 & 0.6 & 5.0   & 112.6  \\
            & 5 & 47.4 & 492.7 & 300.0  & 53.1  & 55.3 & 300.0  & 0.5  & 4.7  & 117.9  \\
          (16/32) & 3 & 50.3 & 488.3 & 300.0  & 58.3  & 50.3  & 300.0  & 1.7 & 4.6  & 249.9  \\
            & 4 & 53.1  & 432.7 & 300.0 & 50.5  & 45.0 & 300.0  & 2.0  & 3.7  & 124.9  \\
            & 5 &  55.1 & 52.0 & 300.0  & 52.4  & 41.0 & 300.0 & 1.2  &  2.3  & 131.3  \\
          (18/36) & 3 & 50.4  & 361.3 & 300.0  & 68.4 &  39.0 &  300.0 & $<$0.1  & 2.2  & 199.2  \\ 
            & 4 & 53.1  & 76.0 & 300.0  &  45.5 & 42.3 & 300.0 &  $<$0.1$^*$  & 2.0$^*$  & 44.8$^*$ \\
            & 5 & 55.2 & 20.0  & 300.0  &  54.3 & 30.3 & 300.0 &  $<$0.1$^*$  & 2.0$^*$ & 52.7$^*$ \\
         (20/40)  & 3 & 54.2$^*$  & 202.5$^*$ & 300.0$^*$  & 79.0  & 38.7 & 300.0  & 0.1$^*$  & 2.0$^*$   & 300.0$^*$  \\ \midrule
         Average & & 37.0 & 874.7&	252.5&	27.2&	62.6&	210.7&	0.3&	8.0&	90.9	  \\ \midrule
    \end{tabular}
\end{table}

As expected, we observe an increase in all metrics with increasing number of nodes and periods. Despite a considerable number of iterations, Benders Decomposition exhibits slow convergence, often struggling to achieve optimality within the 300-minute time limit. 
On average, it reports a 37\% optimality gap and 252.5 minutes solution time. C\&CG performs better, with an average 27.2\% optimality gap and 210.7 minutes solution time, and demonstrates near-optimality within the time limit for more instances. We observe that C\&CG+ is the most favorable approach with an an average of 0.3\% optimality gap within 90.9 minutes solution time. The incorporation of pre-generated columns and constraints at the beginning of the algorithm significantly reduces the number of iterations and enhances algorithm convergence.  Notably, C\&CG+ achieves near-optimality for instances with more than 10 supply nodes, whereas other approaches remain stuck in high optimality gaps. However, from the instance with 18 supply nodes and 4 periods onward, we encounter memory errors in some instances. 
This problem persists across all approaches for every instance involving 20 supply nodes and over 4 periods. These memory errors show the scalability limitations, necessitating exploration of memory-efficient strategies. Despite these challenges, the results clearly demonstrate the superior efficiency and stability of C\&CG+.

\subsection{Case Study}
\label{sec:CaseStudy}
We focus on a case in the Northern Netherlands region where the transition to a green hydrogen economy occurs as part of the project ``Hydrogen Energy Applications in Valley Environments for Northern Netherlands" \citep{Heavenn}. We create our case study based on the expert interviews conducted within the HEAVENN project, recent literature and technical reports for the region. Details of our real-data and instance creation are discussed in Appendix \ref{appendix:CaseStudyDataset}, and 
Figure \ref{fig:map} shows the supply, port and demand nodes on the map of Northern Netherlands based on this dataset.
\begin{figure}[h]
    \centering
    \includegraphics[scale=0.6]{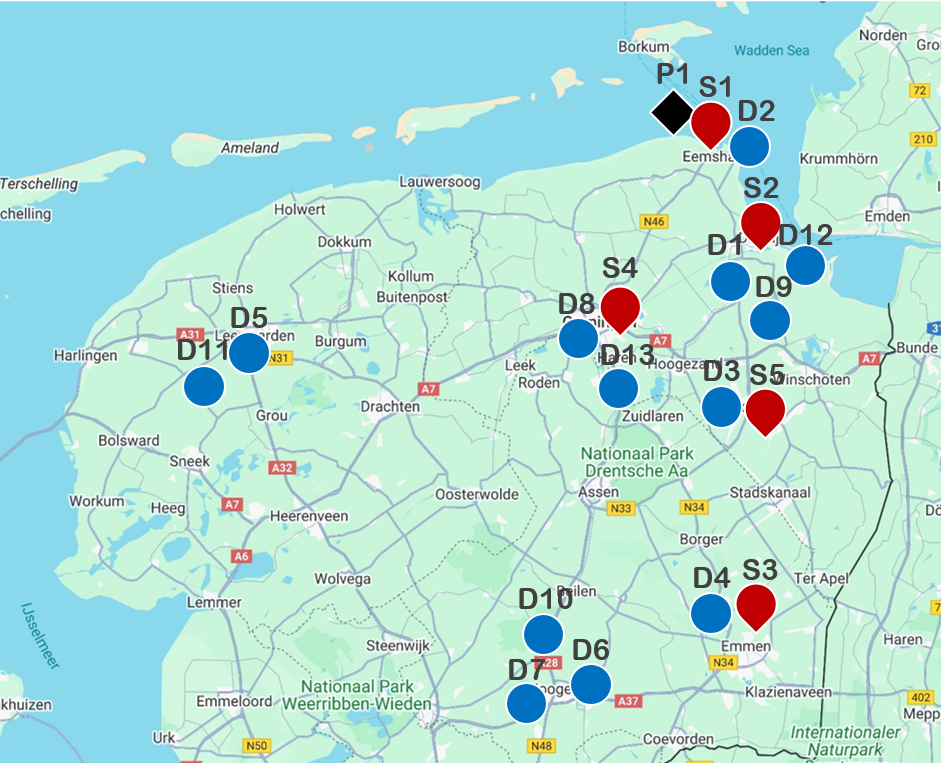}
    \caption{Map of the Northern Netherlands with supply, port and demand nodes marked}
    \label{fig:map}
\end{figure}
We consider planning long-term expansions for the period from 2030 to 2050, making investment decisions every five years – in 2030, 2035, 2040, 2045, and 2050. To create multi-period cost patterns, we rely on available cost projections for the pivotal years, 2030 and 2050. In the absence of direct cost data for the intermediate years, we adopt a pragmatic approach, assuming a linear growth or decline between 2030 and 2050 to estimate costs for these intermediate years. We assume that the cost parameters are fixed within the region across different nodes. In a report by \cite{IEA2019}, the investment cost for per kilowatt of a PEM electrolyzer is estimated to range between €650-€1800 in 2030 and €200-€900 in 2050. Accordingly, we set variable investment cost parameters $c^V_{i1}$ and  $c^V_{i1}$ to the average of the ranges €1225 and €550, respectively. Additionally, we set set-up cost $c^F_{it}$ to €50M. 

We adjust operational cost parameters based on the report by \cite{hyenergy2023}. We set the cost of producing 1 kilogram of hydrogen in 2030 ($c^P_{i1}$) to €3.8 and in 2050 ($c^P_{iT}$) to €2.2. To set transportation cost $c^T_{ijt}$ between nodes $i,j \in \mathcal{N}$ at time $t$, we first calculate the driving distances between nodes using the OSRM API based on real-world map data. Cost of transporting 1 kilogram of hydrogen per km with Gaseous hydrogen tube trailers are expected to drop from €0.01 in 2030 to €0.008 in 2050 \citep{hyenergy2023}. We then multiply the corresponding transportation cost parameters with the calculated distances to obtain $c^T_{ij1}$ for year 2030 and  $c^T_{ijT}$ for year 2050. Previous studies and online sources show significant variations in estimated hydrogen selling prices and import costs for 2030 and 2050. To ensure consistency in our case study, we adjust these parameters in proportion to the production cost. We assume that the selling price of hydrogen $R_{jt}$ is always 20\% higher than the production cost. We set import costs ($c^I_{i1}$) equal to production costs in 2030. However, we anticipate that economies of scale within the HEAVENN region will result in a more substantial decrease in local production costs compared to import costs. Consequently, import costs ($c^I_{iT}$) are projected to be 30\% higher than production costs in 2050. Later in the experiments, we analyze the impact of a different import cost setting on capacity expansion plans.

We set the mean of the hydrogen demand distribution for the years 2030 ($\bar \mu_{j1}$) and 2050 ($\bar \mu_{jT}$) to the values provided for each demand node $j$ in Table~\ref{DemandNodes}. We calculate demand for intermediate years assuming a linear increase during the planning period. We adjust the lower $\underline{\xi}_{jt}$ and upper $\overline{\xi}_{jt}$ bounds of the support sets to $0.75 \times \bar \mu_{j1}$ and $1.25 \times \overline{\xi}_{jt}$, respectively, for each demand node $j \in \mathcal{N}^D$ and time $t \in \mathcal{T}$. 
Expert sessions further emphasize the significance of considering the decision-dependent uncertainty. For the case study, we consider the location-based moment function \eqref{Decision-dependent moment}; however, in Section \ref{Capacity ranges}, we also investigate the effect of the capacity-based moment function \eqref{Decision-dependent capacity moment}. We calculate $\lambda^l_{ij}$ values for each demand node $j \in \mathcal{N}^D$ as a function of distance to supply nodes $i \in \mathcal{N}^S$. Specifically, we use function $\exp\left(-\Delta_{ij}/25\right)$ to assign higher values to shorter distances. We then normalize the obtained values for each demand node by dividing each by the sum of all values. We assume that hydrogen supply openings can increase the demand at a demand node up to 25\%. Therefore, we scale the obtained values considering that their summation $\sum_{i \in \mathcal{N}^S} \lambda^l_{ij} $ is 25\% for each demand node~$j$.

\subsubsection{Value of integrating decision-dependent uncertainty.}
\begin{figure}[h]
    \begin{subfigure}{0.5\textwidth}
        \begin{tikzpicture}[scale=0.8, transform shape]
            \begin{axis}[
                ybar stacked,
                ymax=2000,
                bar width=15pt,
                nodes near coords={},
                enlargelimits=0.15,
                legend style={at={(0.5,-0.20)}, anchor=north, legend columns=-1},
                ylabel={Production capacity (MW)},
                symbolic x coords={2030, 2035, 2040, 2045, 2050},
                xtick=data,
                x tick label style={rotate=45, anchor=east},
                ]
            \addplot+[color4,ybar] coordinates {(2030, 0.0) (2035, 0) (2040, 262) (2045, 272) (2050, 272) };
            \addplot+[color3,ybar] coordinates {(2030, 0.0) (2035, 525) (2040, 632) (2045, 756) (2050, 940) };
            \addplot+[color7,ybar] coordinates {(2030, 0.0) (2035, 0) (2040, 0.0) (2045, 0) (2050, 0) };
            \addplot+[color1, ybar] coordinates {(2030, 0.0) (2035, 0) (2040, 0.0) (2045, 0) (2050, 0) };
            \addplot+[color6,ybar] coordinates {(2030, 0.0) (2035, 152) (2040, 183) (2045, 258) (2050, 326) };
            \legend{\strut S1, \strut S2, \strut S3, \strut S4, \strut S5}
            \end{axis}
        \end{tikzpicture}       
        \captionsetup{font=small}
        \caption{Without DDU}
        \label{fig:WithoutDDUExpansion}
    \end{subfigure}
        \begin{subfigure}{0.4\textwidth}
        \begin{tikzpicture}[scale=0.8, transform shape]
            \begin{axis}[
                ybar stacked,
                bar width=15pt,
                ymax=2000,
                nodes near coords={},
                enlargelimits=0.15,
                legend style={at={(0.5,-0.20)}, anchor=north, legend columns=-1},
                ylabel={Production capacity (MW)},
                symbolic x coords={2030, 2035, 2040, 2045, 2050},
                xtick=data,
                x tick label style={rotate=45, anchor=east},
                ]
            \addplot+[color4,ybar] coordinates {(2030, 0.0) (2035, 0) (2040, 231) (2045, 394) (2050, 394) };
            \addplot+[color3,ybar] coordinates {(2030, 0.0) (2035, 525) (2040, 937) (2045, 1124) (2050, 1341) };
            \addplot+[color7,ybar] coordinates {(2030, 0.0) (2035, 0) (2040, 29) (2045, 38 ) (2050, 48) };
            \addplot+[color1,ybar] coordinates {(2030, 0) (2035, 25) (2040, 39) (2045, 50) (2050, 65) };
            \addplot+[color6,ybar] coordinates {(2030, 113) (2035, 198) (2040, 232) (2045,271) (2050, 316) };          
            \legend{\strut S1, \strut S2, \strut S3, \strut S4, \strut S5}           
            \end{axis}
        \end{tikzpicture}
        \captionsetup{font=small}
        \caption{With DDU}
        \label{fig:withDDUexpansion}
    \end{subfigure}
    \caption{Capacity expansion plans} \label{fig:withvswithoutddu}
\end{figure}
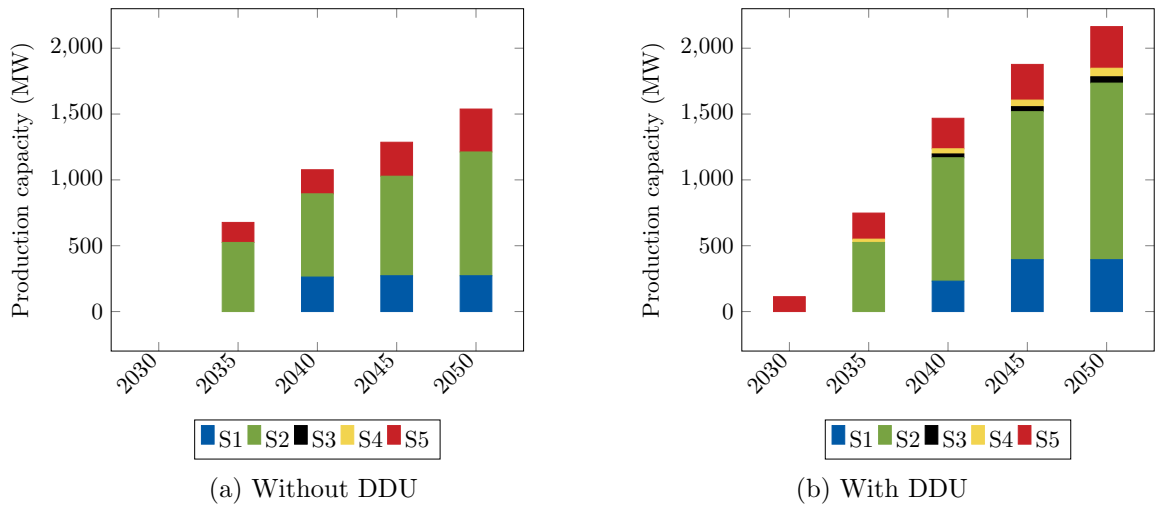

In this section, we analyze the impact of considering decision-dependent uncertainty on capacity expansion plans. We compare two model variants: one with decision-dependent uncertainty (DDU) 
and one model without decision-dependent uncertainty. The model with DDU corresponds to our proposed model. For the model without DDU, we set $\lambda^l_{jt}$ to zero for all demand nodes $j \in \mathcal{N}^D$ at time $t \in \mathcal{T}$. Figure~\ref{fig:withvswithoutddu} illustrates the optimal production capacity decisions for our proposed model with DDU and the model without DDU.

Our findings reveal an increase in total capacity expansion when decision-dependent uncertainty is considered. Specifically, in Figure~\ref{fig:WithoutDDUExpansion}, the model without DDU suggests the construction of 1500MW electrolyzer capacity by the end of the planning period. However, this amount surpasses 2000MW when DDU is considered in Figure~\ref{fig:withDDUexpansion}. Additionally, electrolyzers are built in more nodes when DDU is considered, specifically incorporating Nodes S3 and S4. The DDU model advises higher capacity investments and recommends initiating investments sooner. For example, in Figure~\ref{fig:WithoutDDUExpansion}, we observe hydrogen facilities being built in 2035, while this is already in 2030 in Figure~\ref{fig:withDDUexpansion}. These observations are explained by the anticipated rise in mean demand with the establishment of new electrolyzers. These additional production capacities, including new locations, will meet the increased demand and, consequently, enhance profitability.

\begin{figure}[h]
    \centering
        \begin{tikzpicture}[scale=0.8, transform shape]
            \begin{axis}[
        xlabel={Years},
        ylabel={Production capacity (MW) },
        ymin=0, 
        xtick={2030,2035,2040,2045,2050},
        xticklabel style={/pgf/number format/1000 sep=},
         legend style={at={(1.2,0.5)}, anchor=west},
        ymajorgrids=false,
        grid style=dashed,
    ]

     \addlegendimage{empty legend}
      \addlegendentry{\hspace{-.75cm}\textbf{DDU levels}}
    \addplot[color1,mark=*] coordinates {
 (2030, 0.0)
    (2035, 677)
    (2040, 1078)
    (2045, 1283)
    (2050, 1535)
    };

    \addlegendentry{without DDU}

    \addplot[color2,mark=*] coordinates {
 (2030, 0.0)
    (2035, 677)
    (2040, 1079)
    (2045, 1284)
    (2050, 1535)
    };

    \addlegendentry{10 \% DDU}

            \addplot[color3,mark=*] coordinates {
 (2030, 0.0)
    (2035, 700)
    (2040, 1081)
    (2045, 1568)
    (2050, 1568)
    };

    \addlegendentry{15 \% DDU }
    
    \addplot[color4,mark=*] coordinates {
 (2030, 0.0)
    (2035,702)
    (2040, 1229)
    (2045, 1736)
    (2050, 1805)
    };
    \addlegendentry{20 \% DDU}

            \addplot[color5,mark=*] coordinates {
 (2030, 113)
    (2035, 750)
    (2040, 1471)
    (2045, 1881)
    (2050, 2169)
    };

    \addlegendentry{25 \% DDU}
    \addplot[color6,mark=*] coordinates {
 (2030, 113)
    (2035, 750)
    (2040, 1594)
    (2045, 1894)
    (2050, 2239)
    };
    \addlegendentry{30 \% DDU}

        \addplot[color7,mark=*] coordinates {
 (2030, 113)
    (2035, 1039)
    (2040, 1598)
    (2045, 1897)
    (2050, 2243)
    };
    \addlegendentry{35 \% DDU} 
            \end{axis}
        \end{tikzpicture}
    \caption{Capacity expansions under different levels of decision dependency}\label{fig:DifferentLevels}
\end{figure}
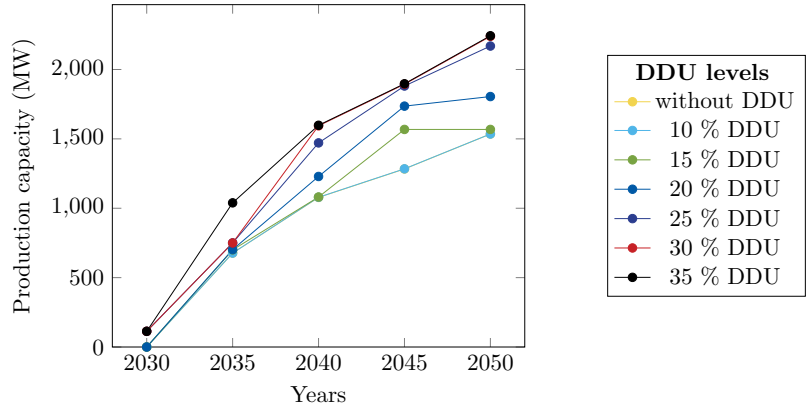

In addition, we analyze the effect of the degree of decision dependency on capacity expansion plans, assuming that the opening of new electrolyzers can increase mean demand by up to a specified level. For the base case, we set this level to 25\%, scaling $\boldsymbol{\lambda^l}$ values so that their summation $\sum_{i \in \mathcal{N}^S} \lambda^l_{ij} $ is 25\% for any demand node $j$. We now vary the decision dependency levels and report the results of the total capacity expansions in Figure~\ref{fig:DifferentLevels}. In the initial stages of planning, capacity expansions are relatively subdued, as the importing cost of hydrogen remains lower than the local production cost. However,  as the decision dependency level increases, we observe an increase in capacity expansions, particularly from 2040 onward. More specifically, we observe that there is no difference between the model without DDU and \%10 DDU level. This indicates that a low level of dependency may not significantly impact decisions.  At 15\% DDU, capacity expansions for the year 2050 shift to 2045. At 20\%, we observe a notable increase in capacity from 2040 onward. As decision dependency levels continue to rise, we observe increasing capacity expansions and shifts to earlier periods. In summary, our results emphasize that varying levels of DDU significantly change optimal capacity expansion plans. The observed shifts in capacity expansions and strategic timelines highlight the critical need for considering decision dependencies. Hence, it is import for businesses to identify the level of decision dependency by analyzing the interaction between investment and demand dynamics to optimally shape their future in the hydrogen sector.

\subsubsection{Out-of-sample Evaluation.}
To evaluate a capacity expansion plan $\hat{\mathbf{x}},\hat{\mathbf{y}}$, we conduct out-of-sample tests by using Evaluation Model~\eqref{model:simulation}. Specifically, we generate a set of demand realizations, based on a given solution $\hat{\mathbf{x}},\hat{\mathbf{y}}$ using the mean information defined in Moment Function~\eqref{Decision-dependent moment}.  We denote these realizations by $d_{jt}^\omega(\hat{\mathbf{x}})$ for each scenario $\omega \in \Omega$, customer site $j \in J$ and time $t \in \mathcal{T}$. For each scenario $\omega$, let $p^\omega $ be the probability of realizing the scenario. Accordingly, we define operational variables $h_{it}^\omega$,  $v_{it}^\omega$ and $z_{ijt}^\omega$  for each scenario $\omega$ representing the amount of hydrogen produced at supply node $i$ at time $t$, the amount of hydrogen imported at port node $i$ at time $t$, and the 
amount of hydrogen transported between nodes $i$, $j$ at time $t$, respectively. Then, we formulate an evaluation model which finds the expected cost over a set of demand realizations under the given investment decisions. The formulation is given as follows:
\begin{subequations} \label{model:simulation}
\begin{align}
\min_{\mathbf{h},\mathbf{z},\mathbf{v}} & \quad \sum_{t \in \mathcal{T}}  \left( \sum_{i \in \mathcal{N}^S} \left( c^{F}_{it} (\hat{x}_{it}- \hat{x}_{i(t-1)} )+ c^{V}_{it} \hat{y}_{it}  \right)  +  \sum_{\omega \in \Omega} p^\omega  \left( \sum_{i \in \mathcal{N}^S} c^{P}_{it}h^\omega_{it} + \sum_{i \in \mathcal{N}^I}  c^{I}_{it} v^\omega_{it}  \right.  \right.  \span \span \nonumber \\
& \hspace{6cm} \left. \left.+ \sum_{j \in \mathcal{N}^D} \left( \sum_{i \in \mathcal{N}^S \cup \mathcal{N}^I} c^{T}_{ijt}   z^\omega_{ijt} -  R_{jt} d_{jt}^\omega(\hat{x}) \right) \right) \right), \span \span \\
\text{s.t.} & \quad   h^\omega_{it}  \leq \sum_{t' \in [1,t]} \hat{y}_{it'},  \hspace{3cm}  &\forall& i \in \mathcal{N}^S, t \in \mathcal{T}, \omega \in \Omega,   \\ 
& \quad h^\omega_{it}  =  \sum_{j \in \mathcal{N}^D} z^\omega_{ijt}, \quad &\forall& i \in \mathcal{N}^S, t \in \mathcal{T}, \omega \in \Omega,    \\ 
&  \quad  v^\omega_{it}  =  \sum_{j \in \mathcal{N}^D} z^\omega_{ijt}, \quad &\forall& i \in \mathcal{N}^I, t \in \mathcal{T}, \omega \in \Omega,   \\ 
& \quad   \sum_{i \in \mathcal{N}^S \cup \mathcal{N}^I  } z^\omega_{ijt} = d_{jt}^\omega(\hat{x}), \quad &\forall& j \in \mathcal{N}^D, t \in \mathcal{T}, \omega \in \Omega,   \\ 
& \quad \mathbf{h}, \mathbf{z}, \mathbf{v} \geq 0.      
\end{align}
\end{subequations}

To evaluate the impact of considering decision-dependent uncertainty and  distributionally robustness, we compare our original setting, referred as \texttt{DRO+DDU}, with three benchmark cases:  \texttt{DRO}, \texttt{DET+DDU}, and \texttt{DET}.  To set \texttt{DRO}, we disregard decision dependency by setting parameter $\lambda^l_{ij}$ to zero for all $i \in \mathcal{N}^S$ and $j \in \mathcal{N}^D$ in Model~\eqref{Continuous Reformulation}. For \texttt{DET+DDU}, we keep decision dependency but exclude distributional uncertainty, constructing a deterministic model explicitly detailed in Appendix~\ref{appendix: deterministic}. Lastly, for \texttt{DET}, we remove decision dependency from the deterministic model. For each capacity expansion plan,  we execute Evaluation Model~\eqref{model:simulation} under 1000 different demand scenarios obtained from a normal distribution whose mean is calculated as $\bar \mu_{jt} \left(1 + \sum_{i \in \mathcal{N}^S} \lambda^{l}_{ij} \hat{x}_{it} \right)$ for each $j \in \mathcal{N}^D$ and $t \in \mathcal{T}$. We present the results in Table~\ref{Tab: simulation}, which includes the average, percentile, and Conditional value-at-risk (CVaR)  values for each solution's out-of-sample objective value, indicating the total cost minus revenue. Lastly, we report the total investment cost including the set-up and variable capacity costs.

Table~\ref{Tab: simulation} shows that \texttt{DRO+DDU} has the best performance in terms of profit. The effectiveness of this approach becomes increasingly apparent with rising CVaR values. \texttt{DET+DDU} stands as the second-best performing approach, underscoring the significance of incorporating the DDU setting. Interestingly, the performance of \texttt{DRO} is comparable to  \texttt{DET+DDU} at CVaR(90\%). This shows the significance of considering both distributionally robustness and DDU when the emphasis is on the worst-case performance evaluation. As expected, the worst approach is \texttt{DET} where both distributionally robustness and DDU are neglected. In deterministic cases, the higher investment costs show a tendency to invest more based on mean demand. The distributionally robust approaches, being more cautious about investments, outperform their deterministic counterparts.
\begin{table}[h]  
\centering
\caption{Out-of-sample results under normal distribution (in millions)} 
\begin{tabular}{p{2cm}p{3.1cm}p{3.1cm}p{3.1cm}p{3.1cm}} 
\toprule
& \texttt{DRO+DDU} & \texttt{DRO} & \texttt{DET+DDU} & \texttt{DET}  \\
\midrule
Average & -2691.13 & -2229.91 & -2452.71 & -2023.73 \\
CVaR(50\%) & -2234.59 & -1863.75 & -1999.74 & -1661.52 \\
CVaR(75\%) & -1984.94 & -1663.29 & -1733.36 & -1445.89 \\
CVaR(90\%) & -1737.44 & -1463.33 & -1464.47 & -1230.79 \\
50\% & -2686.48 & -2222.55 & -2457.39 & -2028.51 \\
75\% & -2292.58 & -1902.79 & -2051.41 & -1704.54 \\
90\% & -1966.63 & -1648.23 & -1732.67 & -1433.54 \\
Investment & 1374.15 & 920.40 & 1571.58 & 1109.19  \\
\midrule
\end{tabular}
\label{Tab: simulation}
\end{table}
\subsubsection{Impact of Import Costs}
We investigate the role of production and import cost dynamics in shaping the capacity expansion plans. When import costs are competitive with local production costs, there is a shift towards imported hydrogen, which reduces investments in local production capacities. However, accessibility to local supply increases potential demand from customers, increasing the overall profit. Therefore, analysing the trade-off between these dynamics is crucial for making investment plans.

\begin{figure}[h]
    \begin{subfigure}{0.5\textwidth}
        \begin{tikzpicture}[scale=0.8, transform shape]
            \begin{axis}[
                ybar stacked,
                ymax=2000,
                bar width=15pt,
                nodes near coords={},
                enlargelimits=0.15,
                legend style={at={(0.5,-0.20)}, anchor=north, legend columns=-1},
                ylabel={Production capacity (MW)},
                symbolic x coords={2030, 2035, 2040, 2045, 2050},
                xtick=data,
                x tick label style={rotate=45, anchor=east},
                ]
                    \addplot+[color4,ybar] coordinates {(2030, 0.0) (2035, 0.0) (2040, 0.0) (2045, 0.0) (2050, 0.0)};
\addplot+[color3,ybar] coordinates {(2030, 0.0) (2035, 0.0) (2040, 0.0) (2045, 665.532) (2050, 665.532)};
\addplot+[color7,ybar] coordinates {(2030, 0.0) (2035, 0.0) (2040, 0.0) (2045, 0.0) (2050, 0.0)};
\addplot+[color1,ybar] coordinates {(2030, 0.0) (2035, 0.0) (2040, 0.0) (2045, 0.0) (2050, 0.0)};
\addplot+[color6,ybar] coordinates {(2030, 0.0) (2035, 0.0) (2040, 155.666) (2045, 186.248) (2050, 221.263)};

            \legend{\strut S1, \strut S2, \strut S3, \strut S4, \strut S5}
            \end{axis}
        \end{tikzpicture}       
        \captionsetup{font=small}
        \caption{Without DDU}
        \label{fig:WithoutDDUImport}
    \end{subfigure}
        \begin{subfigure}{0.4\textwidth}
        \begin{tikzpicture}[scale=0.8, transform shape]
            \begin{axis}[
                ybar stacked,
                bar width=15pt,
                ymax=2000,
                nodes near coords={},
                enlargelimits=0.15,
                legend style={at={(0.5,-0.20)}, anchor=north, legend columns=-1},
                ylabel={Production capacity (MW)},
                symbolic x coords={2030, 2035, 2040, 2045, 2050},
                xtick=data,
                x tick label style={rotate=45, anchor=east},
                ]
\addplot+[color4,ybar] coordinates {(2030, 0.0) (2035, 0.0) (2040, 0.0) (2045, 0.0) (2050, 0.0)};
\addplot+[color3,ybar] coordinates {(2030, 0.0) (2035, 0.0) (2040, 0.0) (2045, 665.671) (2050, 665.671)};
\addplot+[color7,ybar] coordinates {(2030, 0.0) (2035, 0.0) (2040, 0.0) (2045, 0.0) (2050, 0.0)};
\addplot+[color1,ybar] coordinates {(2030, 0) (2035, 46.401) (2040, 60.081) (2045, 70.038) (2050, 70.038)};
\addplot+[color6,ybar] coordinates {(2030, 0) (2035, 121.384) (2040, 143.435) (2045, 247.194) (2050, 247.194)};
            \legend{\strut S1, \strut S2, \strut S3, \strut S4, \strut S5}           
            \end{axis}
        \end{tikzpicture}
        \captionsetup{font=small}
        \caption{With DDU}
        \label{fig:withDDUImport}
    \end{subfigure}
    \caption{Capacity expansion plans with import costs  competitive to production cost} \label{fig:withvswithoutdduimport}
\end{figure}
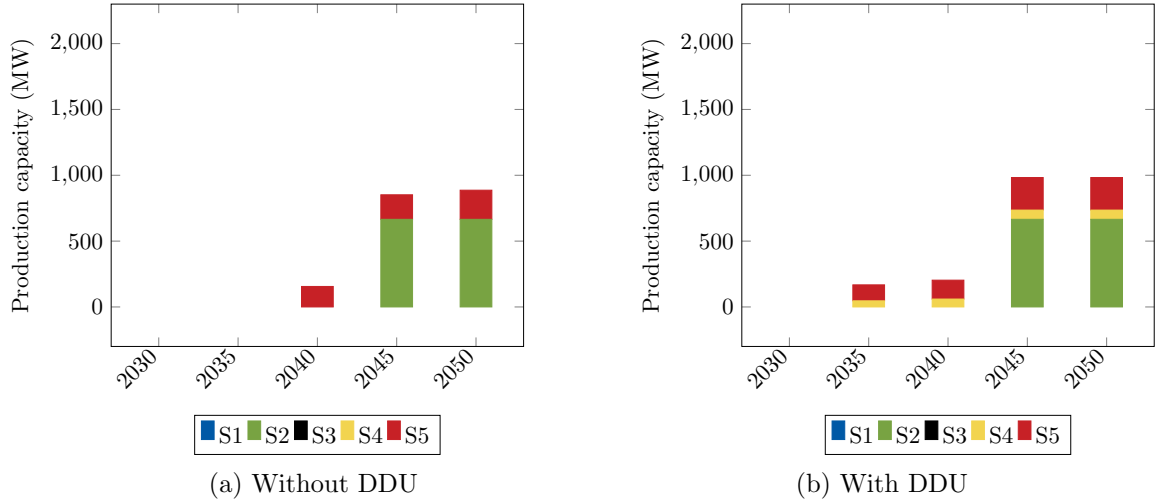

In the base case, import costs are 30\% higher than production costs at the end of the planning period. We adjust this to only 10\% for increasing the competitiveness of imports throughout the planning period.
We report capacity expansion decisions in Figure \ref{fig:withvswithoutdduimport} for without and with DDU cases. Compared to Figure \ref{fig:withvswithoutddu}, we observe decrease in the number of locations and production capacities as expected. Despite these reductions, locations S2 and S5 remain critical due to their proximity to key demand areas D1, D2, and D3. 
Comparing Figures \ref{fig:WithoutDDUImport} and \ref{fig:withDDUImport}, we can say that consideration of DDU still promotes local investments. Competitive import costs may not diminish local investments as significantly as anticipated. It is crucial to thoroughly analyze all related dynamics to ensure that opportunities for local investment are not overlooked.

\subsubsection{Capacity Expansion Plans under Capacity-based Moment Function.} \label{Capacity ranges}
In this section, we analyze the capacity-based reformulation \eqref{Cap:obj} in Appendix \ref{Appendix: capacity-based}. This formulation detects the range in which the production capacity of supply source $i \in \mathcal{N}^S$ at time $t \in \mathcal{T}$ falls and uses this information to affect the demand distribution at demand node $j \in \mathcal{N}^D$. In practice, it requires effort to analyze the effect of production capacity on demand distributions and to determine the corresponding ranges, which may vary based on specific characteristics such as location and time.

\begin{figure}[h]
    \centering
        \begin{tikzpicture}[scale=0.8, transform shape]
            \begin{axis}[
        xlabel={Years},
        ylabel={Production capacity (MW) },
        ymin=0, 
        xtick={2030,2035,2040,2045,2050},
        xticklabel style={/pgf/number format/1000 sep=},
         legend style={at={(1.2,0.5)}, anchor=west},
        ymajorgrids=false,
        grid style=dashed,
    ]
     \addlegendimage{empty legend}
      \addlegendentry{\hspace{-.75cm}\textbf{Ranges for production capacities}}
    \addplot[color2,mark=*] coordinates {(2030, 200.0)
            (2035, 901.1409563776203)
            (2040, 1442.8137448155705)
            (2045, 1854.3979434290416)
            (2050, 2142.1031657090366)};
        \addlegendentry{[0,100],[100,200],[200,$\infty$)}
        \addplot[color3,mark=*] coordinates {(2030, 0.0)
            (2035, 692.7304065547453)
            (2040, 1540.910158119748)
            (2045, 1937.089702360472)
            (2050, 2150.684807623609)};
        \addlegendentry{[0,200],[200,400],[400,$\infty$)}
        \addplot[color5,mark=*] coordinates {(2030, 0)
            (2035, 744.2714361688979)
            (2040, 1181.8255911947886)
            (2045, 1938.5430072631318)
            (2050, 2158.4333782427193)};
        \addlegendentry{[0,300],[300,600],[600,$\infty$)}
        \addplot[color6,mark=*] coordinates {(2030, 0.0)
            (2035, 677.5175177034403)
            (2040, 1098.5939149306989)
            (2045, 1303.606742557932)
            (2050, 1555.1602907374242)};
        \addlegendentry{[0,400],[400,800],[800,$\infty$)}
            \end{axis}
        \end{tikzpicture}
    \caption{Capacity expansions under different ranges for production capacities}\label{fig:DifferentRanges}
\end{figure}
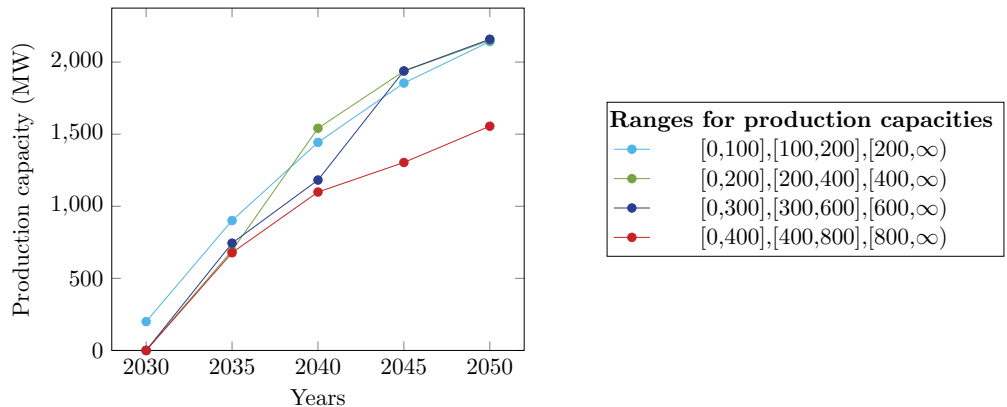

To demonstrate the impact of capacity-based decision dependency in our case study, we make the assumption that there are three ranges of capacities for any supply source $i \in \mathcal{N}^S$ at time $t \in \mathcal{T}$. Similar to the location-based case, we calculate $\lambda^c_{ijr}$ values for each demand node $j \in \mathcal{N}^D$ using the function $\exp\left(-\Delta_{ij}/25\right)$. We then scale the obtained values for each demand node $j$, considering that their summation $\sum_{i \in \mathcal{N}^S} \lambda^c_{ijr}$ is 0\%, 15\%, and 30\% for ranges 1, 2, and 3, respectively. This means that if the capacity of a supply source $i$ at time $t$ falls in the first range, it has no effect on the demand distributions. It starts impacting demand distributions if it falls in the second range and has an even greater impact in the third range. For each range $r \in \{1,2,3\}$, we consider different lower and upper bounds $[l_{itr}, u_{itr}]$ and report the corresponding capacity expansion plans in Figure~\ref{fig:DifferentRanges}. 

First, we observe that the consideration of different capacity ranges significantly alters the capacity expansion plans. In the first case ([0,100],[100,200],[200,$\infty$)), since even lower capacities affect the demand growth, the investments start earlier than in other cases. The second and third cases show similar dynamics. However, the second case makes more investments around the year 2040 compared to the third case. Despite this, by the year 2050, we see that these three cases converge to similar capacity levels. The last case ([0,400],[400,800],[800,$\infty$)) displays a distinct pattern from the others. It resembles a scenario where decision dependency is not considered because it requires higher capacity investments to affect demand. As a result, the system refrains from making further investments. 
For decision-makers, it is crucial to analyze how production capacities influence demand growth in real-life scenarios. This analysis helps identify specific ranges and their effects, thereby enabling more informed decisions regarding future investments.
\section{Conclusions} \label{sec:Conclusions} 
In this paper, we propose a distributionally robust decision-dependent hydrogen network expansion planning problem to deal with uncertainty of hydrogen market and chicken-and-egg dilemma. We propose a model linking hydrogen demand dynamics to investment decisions and offer two reformulations for continuous and discrete support sets. To solve the reformulations, we enhance C\&CG algorithm through solving an inexact master problem, integrating pre-generated columns and constraints, and decomposing the subproblems. 


Computational experiments show the superiority of our solution approach over classical C\&CG and Benders decomposition. Our approach solves over 75\% of instances within the time limit, while Benders Decomposition and C\&CG solve only 35\% and 41\%, respectively. We validate our approach with a case study in the Northern Netherlands within the HEAVENN project by providing insights into the effect of decision dependency on investment plans, out-of-sample performance against benchmark approaches, sensitivity analyses on problem parameters, and analyses on the location-based and capacity-based moment functions. We show that considering the chicken-and-egg dilemma under distributional uncertainty leads to earlier and more diverse investment decisions, making the hydrogen network more profitable.



\bibliographystyle{apalike}
\bibliography{references}

\newpage

\setcounter{page}{1}

\begin{APPENDIX}{}
\section{Reformulation for Bounded Moment Functions}\label{Appendix: Bounded}
We propose a reformulation under the following bounded location-based moment function:
\begin{equation} \label{BoundedMoment}
\mu_{jt}(\mathbf{x})= \min \left( \bar \mu_{jt} \left(1 + \sum_{i \in \mathcal{N}^S} \lambda^{l}_{ij} x_{it} \right), \bar \mu_{jt} \left(1 + \overline{B}_{jt} \right)  \right).
\end{equation} 

Subsequently, we replace the bounded moment function~\eqref{BoundedMoment} in the objective function of Model~\eqref{EquivalentRobustModel} as follows: 
\begin{align} 
\min_{\mathbf{x},\mathbf{y},\alpha,\boldsymbol{\beta}^1, \boldsymbol{\beta}^2,\boldsymbol{\Phi}} & \quad  \sum_{t \in T} \sum_{i \in \mathcal{N}^S} \left( c^{F}_{it} (x_{it}- x_{i(t-1)} )+ c^{V}_{it} y_{it} \right)  + \alpha + \sum_{j \in \mathcal{N}^D} \sum_{t \in \mathcal{T}} \left(  \bar \mu_{jt} \beta^1_{jt} + \bar \mu_{jt} \min\left(\left(\sum_{i \in \mathcal{N}^S} \lambda^{l}_{ij} x_{it}\right) \beta^1_{jt}, \overline{B}_{jt} \beta^1_{jt}\right)   \right. \nonumber \\
& \hspace{3.5cm} \left. + \epsilon_{jt} \beta^1_{jt}  -  \bar \mu_{jt} \beta^2_{jt} - \bar \mu_{jt} \min\left(\left(\sum_{i \in \mathcal{N}^S} \lambda^{l}_{ij} x_{it}\right) \beta^2_{jt}, \overline{B}_{jt} \beta^2_{jt}\right)   + \epsilon_{jt}\beta^2_{jt} \right).
\end{align}

We define $\Upsilon^m_{jt}$ for all $m \in \{1,2\}, j \in \mathcal{N}^D, t \in \mathcal{T}$ to represent $ \min\left(\left(\sum_{i \in \mathcal{N}^S} \lambda^{l}_{ij} x_{it}\right) \beta^m_{jt}, \overline{B}_{jt} \beta^m_{jt}\right)$, and introduce the binary variable $a_{jt}$  for each $j \in \mathcal{N}^D, t \in \mathcal{T}$, which takes the value 1 if $ \sum_{i \in \mathcal{N}^S} \lambda^{l}_{ij} x_{it} $ exceeds $ \overline{B}_{jt}$, and 0 otherwise. We then introduce the following constraints:
\begin{subequations}
\begin{align}
& \sum_{i \in \mathcal{N}^S} \lambda^{l}_{ij} x_{it} \leq \overline{B}_{jt}  + M a_{jt}, & \forall j \in \mathcal{N}^D, t \in \mathcal{T}, \label{Bounded:first} \\
& \sum_{i \in \mathcal{N}^S} \lambda^{l}_{ij} x_{it} \geq \overline{B}_{jt} - M(1-a_{jt}), & \forall j \in \mathcal{N}^D, t \in \mathcal{T}, \label{Bounded:Unnecessary}  \\
& \sum_{i \in \mathcal{N}^S } \lambda^l_{ij}\Phi^1_{ijt}  \leq \Upsilon^1_{jt} + M a_{jt}, & \forall j \in \mathcal{N}^D, t \in \mathcal{T},  \\
& \overline{B}_{jt} \beta^1_{jt} \leq \Upsilon^1_{jt} + M (1-a_{jt}), & \forall j \in \mathcal{N}^D, t \in \mathcal{T},  \label{Bounded:firstlast}    \\
& \sum_{i \in \mathcal{N}^S } \lambda^l_{ij}\Phi^2_{ijt}  \geq \Upsilon^2_{jt}, & \forall j \in \mathcal{N}^D, t \in \mathcal{T}, \label{Bounded:secondfirst}  \\
& \overline{B}_{jt} \beta^2_{jt} \geq \Upsilon^2_{jt}, & \forall j \in \mathcal{N}^D, t \in \mathcal{T}, \label{Bounded:secondlast}  \\
& \Upsilon^1_{jt},\Upsilon^2_{jt} \geq 0, \  a_{jt} \in \{0,1\} & \forall j \in \mathcal{N}^D, t \in \mathcal{T}. \label{Bounded:last}  
\end{align}\end{subequations}

Constraints \eqref{Bounded:first}--\eqref{Bounded:firstlast} represent the positive minimization term in the objective for $m=1$, requiring auxiliary binary variables for their reformulation. In contrast, constraints \eqref{Bounded:secondfirst}--\eqref{Bounded:secondlast} represent  the negative minimization term for $m=2$ and do not need auxiliary binary variables.

\begin{thm}
An equivalent formulation to Model~\eqref{RobustModel} under bounded location-based moment function~\eqref{BoundedMoment}  can be obtained as follows:  
\begin{align}
\min_{\mathbf{x},\mathbf{y},\alpha,\boldsymbol{\beta}^1, \boldsymbol{\beta}^2,\boldsymbol{\Phi}} & \quad  \sum_{t \in T} \sum_{i \in \mathcal{N}^S} \left( c^{F}_{it} (x_{it}- x_{i(t-1)} )+ c^{V}_{it} y_{it} \right)  + \alpha + \sum_{j \in \mathcal{N}^D} \sum_{t \in \mathcal{T}} \left(  \bar \mu_{jt} \beta^1_{jt} + \bar \mu_{jt} \Upsilon^1_{jt} + \epsilon_{jt} \beta^1_{jt} \right. \nonumber \\
& \hspace{9cm} \left.  -  \bar \mu_{jt} \beta^2_{jt} - \bar \mu_{jt} \Upsilon^2_{jt}   + \epsilon_{jt}\beta^2_{jt} \right), \\
\text{s.t.} & \quad \eqref{RobustModel_C1}-\eqref{RobustModel_C3}, \eqref{DI: C2}, \eqref{MC1}-\eqref{MC4}, \eqref{SP:C1}-\eqref{SP:C5}, \eqref{DUALSP:C1}-\eqref{DUALSP:C5}, \eqref{CSBounds},\eqref{Linear_CS1}-\eqref{Linear_CSbound},  \eqref{Continuous Cut}, \eqref{Bounded:first}-\eqref{Bounded:last}.  \nonumber 
\end{align}
\end{thm}

\section{Capacity-based Reformulation}\label{Appendix: capacity-based}
We propose a reformulation under the following capacity-based moment function:
\begin{equation} \label{CapacityBasedMoment}
\mu_{jt}(\mathbf{e})=  \bar \mu_{jt} \left(1 + \sum_{i \in \mathcal{N}^S} \sum_{r \in \mathcal{A}_{it}} 
 \lambda^{c}_{ijr} e_{itr} \right).
\end{equation}

We replace capacity-based moment function~\eqref{CapacityBasedMoment} into Model~\eqref{EquivalentRobustModel} as follows:
\begin{align}
    \min_{\mathbf{x},\mathbf{y},\alpha,\boldsymbol{\beta}^1, \boldsymbol{\beta}^2, \mathbf{e}} & \quad \sum_{t \in \mathcal{T}} \sum_{i \in \mathcal{N}^S} \left( c^{F}_{it} (x_{it}- x_{i(t-1)} )+ c^{V}_{it} y_{it} \right) +  \alpha + \sum_{t \in \mathcal{T}} \sum_{j \in \mathcal{N}^D}  \left( \bar \mu_{jt} \beta^1_{jt} + (\bar \mu_{jt} \sum_{i \in \mathcal{N}^S} \sum_{r \in \mathcal{A}_{it}} 
 \lambda^{c}_{ijr} e_{itr}) \beta^1_{jt} \right. \nonumber \\
 & \hspace{5cm} \left. + \epsilon_{jt} \beta^1_{jt} - \bar \mu_{jt} \beta^2_{jt} - (\bar \mu_{jt} \sum_{i \in \mathcal{N}^S} \sum_{r \in \mathcal{A}_{it}} 
 \lambda^{c}_{ijr} e_{itr}) \beta^2_{jt} + \epsilon_{jt} \beta^2_{jt} \right),  \label{Cap:obj}\\
    \text{s.t.} & \quad \eqref{RobustModel_C1}-\eqref{RobustModel_C3},\eqref{CapDDU: C1}-\eqref{CapDDU: C3},  \eqref{DI: C1}-\eqref{DI: C2}.  \nonumber 
\end{align}

The objective function of Model~\eqref{Cap:obj} contains bilinear terms $e_{itr} \beta^1_{jt}$ and $e_{itr} \beta^2_{jt}$, where their exact linearization is possible due to the binary nature of the variable $e_{itr}$. For linearization purposes, we introduce new variables $\zeta^1_{ijtr}$ and $\zeta^2_{ijtr}$ for all $ \mathcal{A}_{it}, i \in \mathcal{N}^S, j \in \mathcal{N}^D, t \in \mathcal{T}, r \in $. Let $\bar \beta^l_{jt}$ be an upper bound on the variable $\beta^m_{jt}$ where $m \in \{1,2\}$. Accordingly, we obtain the following McCormick constraints: 
\begin{subequations}\label{Cap_McCormicks}
\begin{align}
\zeta^m_{ijtr} &\leq \bar \beta^m_{jt} e_{itr},  & \forall m \in \{1,2\}, j \in \mathcal{N}^D, i \in \mathcal{N}^S, t \in \mathcal{T},\label{Cap_MC1} \\
\zeta^m_{ijtr} &\geq 0, & \forall m \in \{1,2\}, j \in \mathcal{N}^D, i \in \mathcal{N}^S, t \in \mathcal{T} \label{Cap_MC2}, \\
\zeta^m_{ijtr} &\leq \beta^m_{jt}, & \forall m \in \{1,2\}, j \in \mathcal{N}^D, i \in \mathcal{N}^S, t \in \mathcal{T}, \label{Cap_MC3} \\
\zeta^m_{ijtr} &\geq \beta^m_{jt} - \bar \beta^m_{jt} \left(1-e_{itr}\right), \quad & \forall m \in \{1,2\}, j \in \mathcal{N}^D, i \in \mathcal{N}^S, t \in \mathcal{T}. \label{Cap_MC4}
\end{align}
\end{subequations}

\begin{thm}
An equivalent formulation to Model~\eqref{RobustModel} under capacity-based moment function~\eqref{CapacityBasedMoment}  can be obtained as follows:  
\begin{align}
\min_{\mathbf{x},\mathbf{y},\alpha,\boldsymbol{\beta}^1, \boldsymbol{\beta}^2,\mathbf{e},\boldsymbol{\zeta}} & \quad \sum_{t \in \mathcal{T}} \sum_{i \in \mathcal{N}^S} \left( c^{F}_{it} (x_{it}- x_{i(t-1)} )+ c^{V}_{it} y_{it} \right) +  \alpha + \sum_{t \in \mathcal{T}} \sum_{j \in \mathcal{N}^D}  \left( \bar \mu_{jt} \beta^1_{jt} + \bar \mu_{jt} \sum_{i \in \mathcal{N}^S} \sum_{r \in \mathcal{A}_{it}} 
 \lambda^{c}_{ijr} \zeta^1_{ijtr} \right. \nonumber \\
 & \hspace{5cm} \left. + \epsilon_{jt} \beta^1_{jt} - \bar \mu_{jt} \beta^2_{jt} - \bar \mu_{jt} \sum_{i \in \mathcal{N}^S} \sum_{r \in \mathcal{A}_{it}} 
 \lambda^{c}_{ijr} \zeta^2_{ijtr} + \epsilon_{jt} \beta^2_{jt} \right), \label{Cap_Continuous Reformulation}  \\
\text{s.t.} & \quad \eqref{RobustModel_C1}-\eqref{RobustModel_C3}, \eqref{CapDDU: C1}-\eqref{CapDDU: C3}, \eqref{DI: C2}, \eqref{Cap_MC1}-\eqref{Cap_MC4}, \eqref{Continuous Cut}, \eqref{SP:C1}-\eqref{SP:C5}, \eqref{DUALSP:C1}-\eqref{DUALSP:C5}, \eqref{CSBounds} \eqref{Linear_CS1}-\eqref{Linear_CSbound}. \nonumber 
\end{align}
\end{thm}

Note that we propose a reformulation using the capacity-based moment function \eqref{CapacityBasedMoment}. An alternative reformulation for its bounded version, i.e., function \eqref{Decision-dependent capacity moment}, can be obtained by following the steps outlined in Appendix \ref{Appendix: Bounded}.

\section{Dataset Details of Case Study}
\label{appendix:CaseStudyDataset}

In this section, we discuss the potential supply, import, and demand nodes. Northern Netherlands is expected to experience substantial development to host several electrolyzers generating green hydrogen from renewable energy sources. In this region, several locations well-suited for hydrogen production, such as Eemshaven, Delfzijl, Emmen, Groningen, and Zuidwendig, have been identified \citep{laat2022overview, Gasunie2023}. Next to the local hydrogen supply, Eemshaven serves as a crucial gateway for hydrogen import into the region \citep{hyenergy2023, Gasunie2023}. Table~\ref{SupplyData} provides an overview of the considered supply nodes $i \in \mathcal{N}^S$ and port nodes  $i \in \mathcal{N}^I$. 

\begin{table}[h]
\centering \caption{Candidate supply and port nodes}
\label{SupplyData}
\begin{tabular}{p{1.5cm}p{1.8cm}p{2.2cm}p{4cm}}
\toprule
Node ID & Type  & City & Coordinates   \\
  \midrule
 S1 & Supply  & 
  Eemshaven &
    53.44059,   6.82363   \\
S2 & Supply &  
  Delfzijl &
    53.31919,   6.944017   \\
S3 &  Supply & 
  Emmen &
    52.754560,   6.936359  \\
 S4 & Supply & 
Groningen & 
  53.194191,   6.621600  \\
 S5 &  Supply&
Zuidwending &
  53.080842,   6.928074  \\ 
P1 &  Port &  
  Eemshaven &
  53.45162,   6.82991 \\ \midrule 
\end{tabular}
\end{table}

\begin{table}[h]
\centering \caption{Demand nodes}
\label{DemandNodes}
\resizebox{\textwidth}{!}{%
\begin{tabular}{@{}cp{2cm}p{2.2cm}p{5.5cm}p{2.2cm}p{2.2cm}}
\toprule
Node   ID & City  & Coordinates   & Clusters  &  Demand 2030 (ton) & Demand 2050 (ton)   \\ \midrule
D1 & Delfzijl & 53.316718, 6.952066 & \raggedright{Industry cluster in Delzijl} & 120000.00 & 260000.00 \\
D2 & Eemshaven & 53.44026, 6.83999 & \raggedright{Industry cluster in Eemshaven and power plant} & 50000.00 & 90000.00 \\
D3 & Veendam & 53.09345, 6.88291 & \raggedright{Industry cluster in Veendam} & 30000.00 & 60000.00 \\
D4 & Emmen & 52.774836, 6.906929 & \raggedright{Emmtec heating process, GETEC Park} & 150.00 & 400.00 \\
D5 & Leeuwarden & 53.17392, 5.807429 & \raggedright{Hydrogen innovation centre in Leeuwarden} & 4000.00 & 12000.00 \\
D6 & Hoogeveen & 52.72928, 6.507629 & \raggedright{Industry cluster in Hoogeveen} & 3000.00 & 9000.00 \\
D7 & Hoogeveen & 52.707468, 6.411857 & \raggedright{District heating Hoogeveen}& 75.00 & 200.00 \\
D8 & Groningen & 53.206287, 6.598420 & \raggedright{Groningen HRS} & 365.00 & 700.00 \\
D9 & Delfzijl & 53.319264, 6.944039 & \raggedright{Delfzijl HRS} & 32.50 & 55.00 \\
D10 & Hoogeveen & 52.77276, 6.449729 & \raggedright{Hoogeveen HRS} & 32.50 & 55.00 \\
D11 & Leeuwarden & 53.187391, 5.756540 & \raggedright{Friesland HRS}& 32.50 & 55.00 \\
D12 & Delfzijl & 53.321321, 6.942968 & \raggedright{Delfzijl inland ship refueling station} & 14.34 & 20.00 \\
D13 & Groningen & 53.12805, 6.58622 & \raggedright{Groningen airport Eelde} & 16.43 & 27.38 \\
 \midrule 
\end{tabular}}%
\end{table}

Northern Netherlands comprises several key clusters for potential hydrogen demand as detailed in Table~\ref{DemandNodes}. A report by \cite{nibbp} provides insights into the anticipated demand across industry clusters. The report is also supported by HyEnergy Consultancy estimating upwards of 200 kilotons per a year for the entire industry clusters within the region by 2030 \citep{hyenergy2023}. For instance, by 2030, the chemical industry in Delfzijl is expected to require 120 kilotons of hydrogen annually for the production of green methanol and green ammonia. Additional data from \cite{coalition2020northern} supports expected demand for the transportation sector and district heating in Hoogeveen. For example, a pilot project in Hoogeveen aims to transition over 400 homes to hydrogen heating \citep{waterstofwijk2020}. In the context of transportation, previous research has identified potential demand locations including several hydrogen refuelling stations (HRS) and  Groningen airport Eelde \citep{HEAVENN2023}.

\section{Decision-dependent Deterministic Formulation}\label{appendix: deterministic}
\begin{subequations} \label{model:deterministic}
\begin{align}
\min_{\mathbf{x},\mathbf{y}, \mathbf{h}, \mathbf{z}, \mathbf{v}  } & \quad \sum_{t \in \mathcal{T}} \left(\sum_{i \in \mathcal{N}^S} \left( c^{F}_{it} (x_{it}- x_{i(t-1)} )+ c^{V}_{it} y_{it}  + c^{P}_{it}h_{it} \right)  + \sum_{i \in \mathcal{N}^I}  c^{I}_{it} v_{it} \right. \span \span \\ \nonumber
& \hspace{4cm} \left. + \sum_{j \in \mathcal{N}^D} \left( \sum_{i \in \mathcal{N}^S \cup \mathcal{N}^I} c^{T}_{ijt}   z_{ijt} -  R_{jt}   \bar \mu_{jt} \left(1 + \sum_{i \in \mathcal{N}^S} \lambda^{l}_{ij} x_{it} \right)
\right) \right)  \span \span \\
\text{s.t.} & \quad  y_{it}  \leq  M x_{it},  \hspace{4cm} &\forall& i \in \mathcal{N}^S, t \in \mathcal{T},   \\ 
& \quad x_{it}  \geq  x_{i(t-1)},  \quad &\forall& i \in \mathcal{N}^S , t \in \mathcal{T}  \\
& \quad x_{it} \in \{0,1\}, y_{it} \geq 0,  \quad &\forall& i \in \mathcal{N}^S, t \in \mathcal{T}. \\
& \quad   h_{it}  \leq \sum_{t' \in [1,t]} y_{it'},  
&\forall& i \in \mathcal{N}^S, t \in \mathcal{T},   \\ 
& \quad h_{it}  =  \sum_{j \in \mathcal{N}^D} z_{ijt}, \quad &\forall& i \in \mathcal{N}^S, t \in \mathcal{T},    \\ 
&  \quad  v_{it}  =  \sum_{j \in \mathcal{N}^D} z_{ijt}, \quad &\forall& i \in \mathcal{N}^I, t \in \mathcal{T},   \\ 
& \quad   \sum_{i \in \mathcal{N}^S \cup \mathcal{N}^I  } z_{ijt} =  \bar \mu_{jt} \left(1 + \sum_{i \in \mathcal{N}^S} \lambda^{l}_{ij} x_{it} \right), \quad &\forall& j \in \mathcal{N}^D, t \in \mathcal{T},  \\
& \quad \mathbf{h}, \mathbf{z}, \mathbf{v} \geq 0. & 
\end{align}
\end{subequations}

     \end{APPENDIX}

\end{document}